 \documentstyle[10pt,amscd,amssymb]{amsart}

 \def\b{$\bullet$}
 \def\p{\par}

 \theoremstyle{plain}
 \newtheorem{thm}{Theorem}
 \newtheorem{cor}{Corollary}
 \newtheorem{lem}{Lemma}
 \newtheorem{prop}{Proposition}

 \theoremstyle{definition}
 \newtheorem{defn}{Definition}
 
 \newtheorem{rem}{Remark}
 \newtheorem{rems}{Remarks}

 \numberwithin{equation}{subsection}

 \newcommand{\C}{{\cal C}}

 \newcommand{\fg}{{\frak g}}
 \newcommand{\fk}{{\frak K}}
 \newcommand{\fn}{{\frak n}}
 \newcommand{\fb}{{\frak b}}
 \newcommand{\fh}{{\frak h}}

 \newcommand{\ad}{\operatorname{ad}}

 \newcommand{\MOD}{\operatorname{mod}}
 \newcommand{\End}{\operatorname{End}}
 \newcommand{\ann}{\operatorname{Ann}}

 \renewcommand{\hom}{\operatorname{Hom}}

 \newcommand{\Dim}{\operatorname{dim}}
 
 \renewcommand{\ker}{\operatorname{Ker}}

 \newcommand{\Stab}{\operatorname{Stab}}

 \newcommand{\blist}{\begin{list}{\rom{(\roman{enumi})}}{\setlength{\leftmargin}{0em}
 \setlength{\itemindent}{7ex}
 \setlength{\labelsep}{2ex}\setlength{\listparindent}{\parindent}
 \usecounter{enumi}}}
 \newcommand{\elist}{\end{list}}

\begin{document}
 \pagenumbering{arabic}
 \setcounter{page}{1}
\title{Representations of Affine Quantum Function Algebras}
\author{Bharath Narayanan}
\address{Department of Mathematics, University of Arizona}
\maketitle
\section{\bf INTRODUCTION}
\subsection{} The goal of this paper is to investigate the algebraic structure of certain quantized algebras
 of functions associated to affine Kac-Moody Lie algebras and to describe
 their irreducible representations. Let $C$ be an affine Cartan
matrix and $\fg=\fg(C)$ be the associated affine Lie algebra.
 The main object of our interest - ${\Bbb C}_q[G]$ - is a $\star$-subalgebra of the dual space
 $\hom_{k}(U,k)$ generated by matrix coefficients of integrable highest weight $U$-modules, where $\star$ is an involutive
 antiautomorphism.  Similar to the finite-type case, ${\Bbb C}_q[G]$ is a triangular-type algebra whose commutativity relations can be computed using the $R$-matrix.  We construct the irreducible quotient $N_w$ of an induced ${\Bbb C}_q[G]$-module and our main result contains a
description of its annihilator in terms of the Weyl group element
$w$.  Furthermore, these simple modules satisfy a `Tensor-Product
Theorem'
 which asserts that
 $N_w \simeq N_{i_1} \otimes ...\otimes N_{i_k}$, for $w=s_{i_1}...s_{i_k} $
 where the
 $N_i$ are the ${\Bbb C}_{q}[SL_2]$-modules considered as ${\Bbb
 C}_q[G]$-
 modules by restriction to the $i$-th node in the extended Dynkin diagram
 of $g$,
 which induces a surjective homomorphism from ${\Bbb C}_{q}[G]$ to
 ${\Bbb C}_{q}[SL_2]$.  In fact, the modules $N_w \simeq N_{w'}$
  if $w' \simeq w $ in $W$.  Finally, unlike the finite-type case, there is a one-dimensional
   ${\Bbb C}_{q}[G]$-module, $N_{\infty}$, which does not correspond to any Weyl group element.
  The next 3 subsections describe the historical motivations, classical Kac-Moody group theory
  and the finite semisimple quantum function algebras, respectively.
\subsection{} In most approaches to quantum group theory, the basic
object, introduced independently by V.Drinfeld
 \cite {Drin,Drin2} and M.Jimbo \cite {Jimbo}, is the quantized
 enveloping algebra $U_q(\fg)$, which can be viewed as a
 deformation of the universal enveloping algebra $U(\fg)$ of
 a Lie algebra $\fg$.  The quantized function algebra ${\Bbb C}_{q}[G]$ is the {\em non-commutative} or {\em quantum} version of the classical function algebra ${\Bbb C}[G]$. It is a subalgebra of the
dual vector space  $\hom _k(U,k)$, where $k={\Bbb  C}(q)$, $q$ an
indeterminate.  Since the various spectra of the
 coordinate algebra ${\Bbb C}[G]$ of a classical Lie group $G$ contain all the geometric
 information about it, it is natural, from the point of view of
 noncommutative (algebraic) geometry, to study these quantized function algebras.
 A number of techniques have been developed to work with quantized
 function algebras of {\em finite dimensional} semisimple Lie algebras, and there
 are several interesting results and applications. Two important
 non-trivial examples are the discovery of their relation with $q$ special
 functions by L.Vaksman and Y.Soibelman in 1986 and the construction
 of solutions to the Zamolodchikov tetrahedra equations from their
 irreducible representations by D.Kazhdan, M.Kapranov, V.Voevodsky
 and Y.Soibelman in 1992 \cite {KS,KV}.
 \p
However, in the infinite dimensional case, there appears to be
very little literature in this direction \cite {Jo0,Do}, in sharp
contrast to the wealth of information on the representation theory
of $U_q(\fg(C))$ for a Kac-Moody algebra $\fg(C)$ associated to a
Generalized Cartan Matrix $C$ (see \cite {Lu,Jo1,Jo2}). The main
point, even in the finite type case, is that whereas the
irreducible representations of $U(\fg)$ and $U_q(\fg)$ are very
similar, the representation theories for ${\Bbb C}[G]$ and ${\Bbb
C}_{q}[G]$ are extremely different \cite {Yan1}. This difference
stems from the fact that the algebra ${\Bbb C}[G]$ is {\em
commutative}, hence its irreducible representations are one
dimensional, and correspond to the (closed) points of $G$.
\subsection{} It is possible  \cite {KP} to define a Lie group $G=G(C)$ whose Lie
 algebra coincides with the derived subalgebra
 ${\fg}'={\fg}'(C)=[\fg,\fg]$. Denote by $B_{+},B_{-},N_{+},N_{-},H$ and
 $K$ the subgroups
 corresponding to ${\fb}_{+},{\fb}_{-},{\fn}_{+},{\fn}_{-},{\fh}$ and
 ${\fk}$
 respectively.
 In \cite {KP}, the algebra of strongly regular functions ${\Bbb
 C}[G]_{s.r}$
 is
 defined as the algebra generated by matrix coefficients of all
 integrable
 highest weight modules $\overline{L(\Lambda)}$ of the Kac-Moody Lie
 algebra
 ${\fg}(C)$. Their main results are summarized below:
\blist
\item ${\Bbb C}[G]_{s.r.}$ is a unique factorization domain.

\item
 Let $P$ be a subgroup of $G(C)$ and ${\Bbb C}[G]_{s.r.}^{P}$ be the
 algebra of all $f \in {\Bbb C}[G]_{s.r.}$ such that $f(gp)=f(g)$ for
 all
  $p \in P$. Now, let $\theta_{\Lambda}$ be the character of $B_{+}$
  defined by the formula \\$\theta_{\Lambda}$((exp $\thinspace h$)$n$)=$e^{\Lambda (h)}$,
  for $h \in {\fh}$, $n \in N_{+}$ and for any
 $\Lambda \in P_+$, let
  $$S_{\Lambda}=\{ f \in {\Bbb C}[G]_{s.r.} \mid
 f(gb)=\theta_{\Lambda}(b)f(g), \forall g \in G, b \in B_{+} \}.$$ Then,
 \\
  \b(a) (Borel-Weil-type Theorem). The map $\overline{L^{*}(\Lambda)}
 \longrightarrow
  S_{\Lambda}$ defined by $l \mapsto c_{l,v}^{\Lambda}$ is a $G$-
 module isomorphism, where $\overline{L^{*}(\Lambda)}$ is the graded
 dual of $\overline{L(\Lambda)}$.\\
 \b(b) ${\Bbb C}[G]_{s.r.}^{N_{+}}=\oplus_{\Lambda \in P_{+}}S_{\Lambda}$
 and this algebra is isomorphic to ${\oplus}_{\Lambda \in P_{+}}
  \thinspace \overline{L^{*}(\Lambda)}$ as an algebra with the Cartan product:
 $\overline{L^{*}(\Lambda)}$  $\overline {L^{*}(\Lambda')}=\overline{L^{*}(\Lambda + \Lambda')}$, $\forall \Lambda, \Lambda' \in P_{+}$.\\
 \b(c) The algebra ${\Bbb C}[G]_{s.r.}^{N_{+}}$ is a unique factorization
 domain and the coordinate ring of strongly regular functions on
 $\nu_{\Lambda}$ is integrally closed, where $\nu_{\Lambda}$ is
 the (projective) orbit of the highest weight vector $v_{\Lambda}$ in
 $\overline{L(\Lambda)}$ under
  the action of $G$.
 \item (Specmax ${\Bbb C}[G]_{s.r.}) \setminus G$ is non-empty if $C$ is
 of infinite type.
 \item  ${\Bbb C}[G]_{s.r.}$ fails to be a Hopf algebra, since it is neither closed under the comultiplication nor antipodal maps.
\elist

\begin{rems}
\blist
\item
By using only the $U(\fg)$ bimodule structure on $U(\fg)^*$, M.
Kashiwara  \cite {Kashi1}  defines a subalgebra ${\Bbb C}[G]$ of
$U(\fg)^{*}$ satisfying certain finiteness conditions which is
isomorphic to ${\Bbb C}[G]s.r$.  The 'q' analogue of his
definitions are given in Proposition 3.1. 
\item
C. Mokler \cite {Mok} has established that ${\Bbb C}[G]_{s.r.}$ is
really the coordinate ring ${\Bbb C}[M]$ of the monoidal
completion $M$ of $G$, which is isomorphic to ${\Bbb C}[G]$ by
restriction, since $M$ is a monoid containing $G$ as its group of
units. He also computes Specmax ${\Bbb C}[G]_{s.r.}) \setminus G$.
\item The books by O. Mathieu \cite {Mat} and S. Kumar \cite {Kumar} examine the structure of Kac-Moody groups in detail.
\elist
\end{rems}
\subsection{} When $G$ is a finite-dimensional semisimple Lie algebra, ${\Bbb C}_{q}[G]$ is the Hopf-dual of $U_{q}(\fg)$, equivalently defined as the algebra spanned by matrix coefficients of all finite-dimensional $U_{q}(\fg)$-modules (Peter-Weyl type theorem).  It has a triangular structure and its commutativity relations are obtained using the quantum $R$-matrix.
  The main results in the representation theory of ${\Bbb C}_{q}[G]$ \cite {Yan1, Jo1} are summarized below:
\blist
\item The irreducible representations $V_w$ of ${\Bbb C}_{q}[G]$ are parameterized
 by Weyl group elements $w$ and are independent of the reduced presentation of
$w$.  The representations corresponding
 to the identity element in the Weyl group are 1-dimensional, while the
 others are infinite-dimensional \cite {Yan2}.
\item The collection of $V_w$'s is in 1-1 correspondence with the symplectic leaves of
 the classical Poisson-Lie group $G$ associated with $\fg$.
\item
 $V_w$ is the irreducible quotient of the ${\Bbb C}_{q}[G]$ module induced from a one-dimensional ${\Bbb C}_{q}[G/N]$-module, where $N$ is the nilpotent subgroup of $G$
\cite {Jo1}.
\elist

\subsection{} This paper is adapted from my Ph.D thesis \cite {ambi} and I would like to thank Y.Soibelman for suggesting the main result to me.  I also acknowledge here that the construction of the induced modules and their annihilators follows the approach of A.Joseph in \cite {Jo1} for finite
 dimensional semisimple Lie algebras.  Finally, a special word of thanks to my advisors, Z.Lin and A.Rosenberg for their invaluable
 guidance.

 \par
 \vfill\pagebreak

 \section{\bf PRELIMINARIES}
 \subsection{Notations and Basic Definitions}
 Here I summarize the results and notation relating to Kac-Moody algebras
and quantum groups. A more detailed treatment can be found in the books
by Kac \cite {Kac} and Lusztig \cite {Lu}.
 Let $q \in {\Bbb C} \setminus 0$ which is not a root of unity and
 let $k={\Bbb Q}(q) \subseteq {\Bbb C}(q)$.
The following concise description of a Kac-Moody Lie algebra is
adapted from \cite {Chari}.
 \begin{defn}
 Let $C= ((a_{ij}))$ be a symmetrizable $l{\times}l$ generalized Cartan
 matrix of rank $r$. It is defined by  integers $a_{ij}$, with $a_{ij}$
 being non-positive for $i \neq j$
 such that $$a_{ii}=2, a_{ij} = 0 \Rightarrow a_{ji} =0 $$
 and let $d_i$ be the coprime positive integers such that the matrix
 $((d_ia_{ij}))$ is symmetric.  Fix the index set $I=\{1,2,...,l\}$.
 Denote by ${\fg}'(C)$ the complex Lie algebra on $3l$ generators $e_i,f_i,h_i, i \in I$ and the defining relations:
 $$[h_i,h_j]=0, [e_i,f_j]=\delta_{i,j}h_{j},$$
 $$[h_i,e_j]=a_{ij}e_j, [h_i,f_j]=-a_{ij}f_{j},$$
 $$(\ad e_i)^{1-a_{ij}}e_j=0, (\ad f_i)^{1-a_{ij}}f_{j}=0,    i \neq j.$$
 Let ${\fh}'$
be the linear span of the $h_i, i \in I$.  Choose a vector space ${\fh}''$ of
dimension $l-r$, with basis $\{ D_{r+1},...,D_{l} \}$.
 The associated {\em Kac-Moody algebra}, $\fg(C)$ is the Lie algebra with generators $e_i,f_i,h_i,i \in I$ and $D_i, i=r+1,...,l$
and with defining relations those of ${\fg}'(C)$ together with
$$[D_i,D_j]=0, \hspace{.3 cm}  [D_i,h_j]=0, \hspace{.3 cm}   [D_i,e_j]={\delta}_{i,j}e_j, \hspace{.3 cm} [D_i,f_j]=-{\delta}_{i,j}f_j.$$
\end{defn}

\begin{rem}
The direct sum $\fh={\fh}'\oplus {\fh}''$ is called the {\em Cartan subalgebra}
of ${\fg}(C)$.
\end{rem}
 Let $q_i=q^{d_i}$ and define the simple roots  $\alpha_{i}:\fh \longrightarrow {\Bbb C}, i \in I,$ by $$\alpha_{i}(h_{j})=a_{ji},\alpha_{i}(D_{j})=\delta_{i,j}.$$
 They are linearly independent. Let $\pi(C)=\{ \alpha_1,\alpha_2,...\alpha_{l}\}$ denote the set of simple roots of $\fg$.  There is a  non-degenerate bilinear form on ${\fh}^{*}$ such that $d_{i}a_{ij}=(\alpha_i,\alpha_j)$, $i,j \in I$. Write the
 fundamental weights as $\{ \omega _{\alpha} \in {\fh}^{*};\alpha \in
 \pi(C) \}$,
 satisfying $ (\omega_{\alpha},\beta^{\vee})=\delta_{\alpha, \beta} $,
 where $\beta \in {\fh}^{*}$, $\beta^{\vee}:=2\beta/(\beta,\beta)$.
 Let $P(C)=\oplus_{\alpha \in \pi(C)}{\Bbb Z} \omega _{\alpha} + {\Bbb Z}
 \pi(C) \subseteq \{ \lambda \in {\fh}^{*} \mid
 (\lambda,\alpha^{\vee}) \in {\Bbb Z}, \forall \alpha \in \pi(C) \}$ and $P_{+}=\{ \lambda \in
 P(C) \mid (\lambda,\alpha^{\vee}) \in {\Bbb N} \}.$
\p
Define linear maps $s_i:{\fh}^{*} \longrightarrow {\fh}^{*}$, called the simple reflections, by $$s_{i}(\alpha)={\alpha}-{\alpha}(h_{i}){\alpha}_{i}.$$
The Weyl group $W$ of $\fg$ is the subgroup of $GL({\fh}^{*})$ generated by $s_1,...,s_n.$  The action of $W$ preserves the bilinear form $(,)$ on ${\fh}^{*}$.
\p
 Introduce $\rho \in h^{*}$ by
 $(\rho,\alpha_{i}^{\vee})=a_{ii}/2$,  $i \in I$.  Note that $\rho$ is not
unique if det $C=0$ and we pick any
 solution.
 Let $e$ denote the identity element in the Weyl group $W$ and
 $P_{++}=\rho + P_{+}=\{\Lambda \in P_{+} \mid
 {\Stab}_{W}\Lambda=\{e\}\}$.\\
 Define $[n]_{q}=(q^{n}-q^{-n})/(q-q^{-1})$,
 $[n]_{q}!=[n]_{q}[n-1]_{q}...[1]_{q}$, $[0]_{q}!=1$
 and ${n \choose t}_{q}=[n]_{q}!/([t]_{q}![n-t]_{q}!).$
 Denote by $\Phi$ the root
 system of ${\fg}$ with respect to $\fh$.
 Let $$\Gamma=\{ \lambda \in {\Bbb Z}{\Phi} \otimes_{\Bbb Z} {\Bbb Q} \mid
 (\lambda,\mu) \in {\Bbb Z},  \forall \mu \in \sum_{\beta \in \pi (C)} {\Bbb Z} \omega_{\beta} \},$$ where ${\Bbb Z} \Phi = \sum_{\alpha \in \Phi}{\Bbb
 Z}{\alpha}$.
 \begin{defn}
 Let $U=U_q(\fg)$ be the Hopf algebra over $k$
 with generators \\
 $\left<E_{\alpha},F_{\alpha},K_{\lambda}\right> $ ($\alpha \in \pi(C),\lambda \in
 \hom_Z(\Gamma,Z)$) and defining relations:
 \[K_{\lambda}K_{\mu}=K_{{\lambda}+{\mu}},\]
 $$K_{\lambda}E_{\beta}K_{\lambda}^{-1}=q^{(\lambda,\beta)}E_{\beta},$$
 $$K_{\lambda}F_{\beta}K_{\lambda}^{-1}=q^{-(\lambda,\beta)}F_{\beta},$$
 along with the quantum Serre relations: for $i \neq j$,
 \[\sum_{k=0}^{1-a_{ij}}(-1)^{k}{{1-a_{ij}}\choose
 k}_{q_{i}}(E_i)^{k}(E_j)(E_i)^{1-a_{ij}-k}=0.\]
 \[\sum_{k=0}^{1-a_{ij}}(-1)^{k}{{1-a_{ij}}\choose
 k}_{q_{i}}(F_i)^{k}(F_j)(F_i)^{1-a_{ij}-k}=0.\]
 \end{defn}

 Let $i \in I$.  We can identify $\alpha_{i} \in {\hom}_{\Bbb Z}(\Gamma,{\Bbb Z})$
 by $\alpha_{i}(\beta)=(\beta,\alpha_{i}^{\vee})$. Denoting $E_{\alpha_{i}},F_{\alpha_{i}},K_{\alpha_{i}}$
 by
 $E_i,F_i,K_i$ respectively, the Hopf algebra structure on $U$ is given
 by the
 antipode $S$, comultiplication $\Delta$, and counit $\epsilon$, which are
 defined
 on the generators via:
\blist
 \item $ S(E_i) = -K_{i}^{-1}E_{i}$, $ S(F_i)=-F_{i}K_{i}$, $ S(K_i)=K_{i}^{-1}$
 \item $\Delta (E_i)=E_i \otimes 1 + K_i \otimes E_i$, $\Delta (F_i)=F_i \otimes K_i^{-1} + 1 \otimes F_i$, $\Delta (K_i)=K_i \otimes K_i.$
 \item $\epsilon (E_i) = \epsilon (F_i) = 0$,  $\epsilon (K_i) = \epsilon (K_i^{-1}) =1.$
\elist
\begin{rem}
 $U$ is called the {\em simply connected form} of the quantum enveloping algebra.  For each $i \in I$, the subalgebra $U_{i}$ of $U$
 generated by the $E_{i},F_{i}$ and $K_{i}$ is isomorphic to
 $U_{q}(sl_{2})$ as a Hopf algebra.  Denote by $U^{\geq 0}$ (resp. $U^{> 0}$) the subalgebra
 of $U$ generated by the $E_{i}$ and $K_{i}$, $i \in I$ (resp. $E_{i}$,
 $i \in I$).  Similarly define $U^{\leq 0}$ and $U^{<0}$ by replacing
 $E_{i}$ with $F_{i}$.
 \end{rem}

 \subsection{$\star$-structures}
 Classically the $\star$-structure on the Lie algebra $\fg$, with generators
 $\{ e_i,f_i,h_i \}$, is given by a Lie algebra antiautomorphism
 interchanging $e_i$
 with $f_i$ and preserving the $h_i$. Then the (real) compact form of
 $\fg$ is
 $\fk=\{x \in g \mid x^{\star}=-x \}$. Recall that this is the Lie algebra
 of the compact subgroup $K<G$, where $G$ is the Lie group of the Lie algebra
 $\fg$.  Let $U^{0}$ denote the group algebra of the multiplicative group $T$ in $U$
 generated by $\{K_{\alpha}:\alpha \in \pi(C)\}$.
 \begin{defn} Let $ \omega$ be the Cartan involution on $U$. It is an
 algebra
 automorphism given
 by  $$\omega (E_i)=-F_i, \hspace{.3 cm}  \omega(F_i)=-E_i, \hspace{.3 cm}  \omega (K_i) = K_i^{-1}, \hspace{.3 cm}  \omega(q)=q.$$
 \end{defn}
\p
\begin{defn}
Let $A$ be any complex associative algebra.  A $\star$ structure
on $A$ is an antilinear, involutive, algebra antiautomorphism and
coalgebra automorphism:
 $$(cx)^{\star}=\bar c(x^\star), \hspace {0.3 cm} (x^\star)^\star=x, \hspace {0.3 cm} (xy)^{\star}=y^{\star}x^{\star}$$
 for any $c \in {\Bbb C}$ , $x,y \in A$, where $\bar {x}$ is the complex conjugate of $x$.
\end{defn}
$U$ can be given a $\star$-structure by setting $\star=S \circ
\omega.$   The dual vector space $U^{*}$ (with a different meaning
of $^{*}!$) then becomes a $\star$-algebra via
 $l^{\star}(u)=\overline {l(S(u)^{\star})}$, i.e. under the transpose of the
 automorphism $\star \circ S=\omega$ on $U$ where $l \in U^{*} $ and $u \in
 U$.
 The following (equivalent) identities are easily checked :
  \[S(u)^{\star}=S^{-1}(u^{\star}) , S \circ \omega \circ S = \omega. \]

 \subsection{Integrable U-modules and their duals}
  By  $U$-modules we will mean {\em left} $U$-modules, unless otherwise
 specified. For any $U$-module $M$ over $k$, define a
 $U$-module structure (resp. {\em right} $U$-module structure) on the
 dual
 vector space $\hom_k(M,k)$ via
 $(uf)(m)=f(\gamma(u)m)$ (resp.) $(fu)(m)=f(um)$ where $\gamma$ is
 an antiendomorphism of $U$, $u \in U$ and $m \in M$. If we
 take $\gamma = S$, the antipode of $U$, the resulting module is denoted
 $M^{*}$ and if we take $\gamma = \omega \circ S$, we denote it
 $M^{\circ}$.  Further, if $M'$ is any right $U$-module, then
$\hom_k(M',k)$ is a left $U$-module, via
 $(uf)(m)=f(mu)$.
\p
  Now fix $\phi:U \longrightarrow {\End}(M),$
 $M=L(\Lambda)$, highest weight, integrable $U$-module of highest weight $\Lambda \in P_{+}$, with highest weight vector $v_{\Lambda}$.
 Let $L(\Lambda)_{\lambda}^{*}=(L(\Lambda)_{-\lambda})^{*}$ and
$L^{*}(\Lambda) = \bigoplus
_{\lambda \in \Omega (\Lambda)} L(\Lambda)_{\lambda}^{*}.$
$L^{*}(\Lambda)$ is a lowest weight integrable module with lowest weight
vector $l_{-\Lambda}$
 \begin{defn}
 If we give the vector space $L(\Lambda)$ a left $U$-module structure by
replacing $\phi$ by $\phi ^{\sharp} = \phi \circ \omega$, it becomes the
 {\em graded dual}, denoted
 $L(\Lambda)^{\sharp}$.
\end{defn}
\begin{rems}
\blist
\item As $U$-modules, $L(\Lambda)^{\sharp} \simeq L^{*}(\Lambda)$, since
the map taking $v_{\Lambda}$ to $l_{-\Lambda}$ extends to a $U$-module
isomorphism from $L(\Lambda)^{\sharp}$ to $L^{*}(\Lambda)$.
\item $L(\Lambda)^{\sharp} = L(\Lambda)^{*}$ if and only if $L(\Lambda)$ is finite dimensional.
\elist
\end{rems}

\begin{defn}
 Let $M_{1}$ and $M_{2}$ be $U$-modules. A bilinear form $\left<,\right>:M_{1} {\otimes} M_{2} \longrightarrow k$ is said to be $U$-invariant if $\left<uv_1,v_2\right>=\left<v_1,S(u)v_2\right> \forall u \in U$ and $v_i \in M_i, i=1,2$.
 \end{defn}

 There exists a unique $U$-invariant bilinear form on
 $L^{*}(\Lambda) \otimes L(\Lambda)$
 satisfying $$\left<l_{-\Lambda},v_{\Lambda}\right>=1,$$ where $l_{-\Lambda}$ is a
 fixed lowest weight vector for $L(\Lambda)^{\sharp}$.  It is
nondegenerate. Using
 $\left<,\right>$ we may regard $L(\Lambda)^{\sharp}$ as a subspace of the dual
vector space
 $L(\Lambda)^{*}$. Note that if a statement holds for $L(\Lambda)$ (resp.$
L(\Lambda)^{*}$)
then a
 similar one holds for $L(\Lambda)^{\sharp}$ (resp $L(\Lambda)^{\circ}$)
using
 $\omega$(resp $\omega \circ S$).
Finally, observe that if $C$ is of finite type, then $L(\Lambda)$
 is finite dimensional as a vector space and it follows that
 $L(\Lambda)^{\circ} \simeq L(-w_{0}\Lambda)$ and $L(\Lambda)^{\sharp}
 \simeq L(-w_{0}\Lambda)^{*}$, as $U$-modules, where $w_{0}$ is the Weyl
 group element of maximal length.
\par
\vfill\pagebreak

 \section{\bf QUANTIZED FUNCTION ALGEBRAS OF A KAC-MOODY LIE ALGEBRA}

 \subsection{Definition of $A={\Bbb C}_{q}[G]$}
Since $U$ is a coalgebra, its vector space dual $U^{*}$ is naturally an
algebra.
 \begin{defn} The {\em matrix coefficient} $C_{l,v}^{\Lambda}$ of an
 integrable
 highest weight (type 1) $U_q(\fg)$-module $L(\Lambda)$ is defined via:
 $$C_{l,v}^{\Lambda}(u):=\left<l,\rho_{\Lambda}(u)v\right>,$$
 for $l \in L(\Lambda)^{\sharp}$, $v \in L(\Lambda)$ and $u \in
 U_q{\fg}.$
 \end{defn}
  Let $F:=\bigoplus_{\Lambda \in P_{+}}L(\Lambda)^{\sharp}
 \otimes L(\Lambda)$. Then the map $l \otimes v \mapsto
 C_{l,v}^{\Lambda}$ extends to a vector space map $F \longrightarrow
 {\hom}_{k}(U,k).$ Let $R$ denote its image. It carries a multiplicative
 structure
 via:
 $$(C_{l,v}^{\Lambda}C_{l',v'}^{\Lambda'})(u)=(l \otimes
 l')((\phi_{\Lambda} \otimes {\phi_{\Lambda'}})\Delta(u) (v \otimes v')).$$
 This map is well defined since the tensor product of two integrable highest weight modules decomposes as the direct sum of integrable highest weight modules.

 \begin{defn} Let $R$ be the subspace of $U^{*}=\hom_{k}(U,k)$
 spanned by matrix coefficients $C_{l,v}^{\Lambda}$ with $\Lambda \in
 P_{+}$.
 \end{defn}
 \begin{rem}This is the quantum analog of the algebra of strongly regular functions on the group $G$ associated to the  derived subalgebra ${\fg}' =
 [{\fg},{\fg}]$, $G$ being
 an infinite-dimensional affine algebraic group of Shafaravich type
 \cite {KP}. It is important to note that, {\em unless $C$ is of finite type},
 $R$
 is {\em not a Hopf algebra}, since it is {\em not} closed
 under the comultiplication (resp. antipodal) map dual to the
 multiplication (resp. antipodal) map in $U$. The {\em Hopf Dual} of
 $U$, denoted {\bf $U^{\star}$}, is defined \cite {Jo1} as $U^{\bullet}:=$\{$\xi \in U^{*} \mid \xi (I)=0$  for some
 two-sided ideal $I$ in $U$ of finite codimension\}. It is a fact
that $U^{\bullet}$ consists of matrix coefficients of finite
dimensional $U$-modules and is properly contained in $R$ in the
infinite type case.
 $U^{\bullet}$ coincides with $R$ if $C$ is of finite type and is commonly referred to in the literature as the {\em Quantum
 Coordinate Algebra} of $G$ \cite {APW}.
 \end{rem}
\begin{defn}
Letting $U_{i}^{\geq 0}:=U^{\geq 0}U_{i}$ denote the quantum parabolic algebra, define
\blist
\item $\cal C$=\{$J \mid J$ is a 2-sided ideal of $U_{i}^{\geq 0}$, $\dim (U_{i}^{\geq 0} / J) < \infty$, and $U^{0}$ acts semisimply on $(U_{i}^{\geq 0} / J)$ with weights in $P(\pi), \forall i \in I$\},
\item $\cal F$ = $\oplus_{\Lambda \in P_{+}}$\{$\phi \in \hom _{k}(U,k) \mid \phi$ is $U_{i}^{\geq 0}$-finite under left and right multiplication for all $i \in I$ and  $\phi$ is a weight vector of weight $\Lambda$ under the left action of $U^{0}$\}.
\elist
\end{defn}
The proof of the next proposition is straightforward.
\begin{prop}
Let $\phi \in \hom _{k}(U,k)$. Then the following statements are equivalent:
\blist
\item $\phi \in R$
\item $\phi(J)=0$ for some $J \in \cal C$.
\item $\phi \in \cal F$.
\elist
\end{prop}
\p Recall (Section 2.2) that $\hom _{k}(U,k)$ is a
$\star$-algebra, with the involution $\star$ defined via
$l^{\star}(u)=l(S(\omega(S(u))))=l(\omega(u))$, for any $l \in
\hom _{k}(U,k)$, $u \in U$.

\begin{defn} Define the {\em quantized function algebra of a Kac-Moody Lie algebra}, ${\Bbb C}_{q}[G]:=\left<R,\star(R)\right>$ as  the minimal $\star$-subalgebra of
 $Hom_{k}(U,k)$ generated by elements of $R$, where $\star$ is the
 transpose of the involution $\omega$ on $U$, and denote it by $A$.
 \end{defn}

 \begin{rem} If we interchange $U^{>0}$ and $U^{<0}$ in the definition of
highest weight module to define a {\em lowest} weight module,
it follows from Proposition 3.2 that the image of a
 matrix coefficient of a {\em highest} weight $U$-module, under $\star$, is
 a matrix coefficient of a {\em lowest} weight $U$-module.  Thus $A={\Bbb C}_{q}[G]$ can be viewed as the analogue of the algebra of holomorphic and anti-holomorphic functions of a complex variable.
 \end{rem}

 \subsection{U-bimodule structure}
 The elements $C_{l,v}^{\Lambda}$ of $R$ are written as
 $C_{-\lambda,i;\mu,j}^{\Lambda}$ if
 $l \in L(\Lambda)^{\sharp}_{\lambda}$ is the $i$th basis vector $1 \leq
 i \leq {\dim}(L(\Lambda)_{\lambda}) $ and $v \in L(\Lambda)_{\mu}$ is the
 $j$th basis vector $1 \leq j \leq \dim(L(\Lambda)_{\mu})$. It is a
 common practice to omit the indices $i$ and $j$, since all formulae
 involved
 are expected to hold irrespective of their choice.

 \begin{prop} Let $u \in U$.
 The following relations express the $U$-bimodule structure of $A$.
 $$uC_{l,v}^{\Lambda}=C_{l,uv}^{\Lambda}, \hspace{.9 cm} C_{l,v}^{\Lambda}u=C_{lu,v}^{\Lambda}, \hspace{.9 cm} u(C_{l,v}^{\Lambda})^{\star}=(C_{l,\omega(u)v}^{\Lambda})^{\star}, \hspace{.9 cm}(C_{l,v}^{\Lambda})^{\star}u=(C_{l\omega(u),v}^{\Lambda})^{\star}.$$
 \end{prop}
 \begin{pf} Use the definition of the $U$-bimodule structure on
 $\hom_{k}(U,k)$ and the right $U$-module structure on
 $\hom_{k}(V,k)$, for any left $U$-module $V$. For example, for any $x
 \in U$,
 we have that
 $$uC_{l,v}^{\Lambda}(x)=C_{l,v}^{\Lambda}(xu)=l(\rho_{\Lambda}(xu)v)=C_{l,uv}^{\Lambda}(x),$$ and
 $$(C_{l,v}^{\Lambda}u)(x)=C_{l,v}^{\Lambda}(ux)=l(\rho_{\Lambda}(ux)v)=lu(\rho_{\Lambda}(x)v)=C_{lu,v}^{\Lambda}(x).$$
 \end{pf}

 \subsection{Triangular-Type Structure}
 \begin{defn} Let $A_{+}$ denote the subspace of $A$ spanned by the $\{C_{-\mu,j;\Lambda}^{\Lambda} ; \Lambda \in
 P_{+},\mu \in \Omega(\Lambda)\}$, and $A_{-}=\star(A_{+})$.  Similarly define $A_{++}$ (resp. $A_{--}$) by replacing $P_{+}$ by $P_{++}$ in the definition of $A_{+}$ (resp. $A_{-}$).

 \end{defn}

 \begin{rem}
We define the set of left $U^{>0}$-invariants of $A$ as the set \{ $\xi \in A \mid x{\xi}=\epsilon(x)\xi, \forall x \in U^{>0}$ \}.  Similarly define the left $U^{<0}$-invariants as \{ $\xi \in A \mid y{\xi}=\epsilon(y)\xi, \forall y \in U^{<0}$ \}.  Then it is clear that $A_{+}$ (resp. $A_{-}$) is just the set of invariants of $A$ with respect to
 the left action of $U^{>0}$ (resp. $U^{<0}$) and may be viewed as the
 quantum function algebras on $G/N_{+}$ (resp. $G/N_{-}$).
 \end{rem}

 \begin{defn}
 For $\Lambda \in P_{+}$, define the subspaces $L^{+}(\Lambda)^{*}$
 (resp. $L^{-}(\Lambda)^{*}$)  of $A_{+}$
 (resp. $A_{-}$) as $$L^{+}(\Lambda)^{*}=\{
 C_{{\xi},{\Lambda}}^{\Lambda} \slash \xi \in L(\Lambda)^{\sharp} \}.$$
 $$L^{-}(\Lambda)^{*}=\{ (C_{{\xi},{\Lambda}}^{\Lambda})^{\star} \slash \xi
 \in
 L(\Lambda)^{\sharp} \}.$$
 \end{defn}

 \begin{rem}
 As a right $U$-module, $A_{+}$ is a direct sum of the $L^{+}(\Lambda)^{*}$ which is isomorphic to $L(\Lambda)^{\sharp}$. Similarly,  $A_{-}$ is a direct sum of the $L^{+}(\Lambda)^{*}$ which is isomorphic to $L(\Lambda)^{\sharp}$. Further, these modules satisfy the Cartan multiplication rule:
 $$L^{+}(\Lambda)^{*}L^{+}(\Lambda')^{*}=L^{+}(\Lambda + \Lambda')^{*}, \forall \Lambda, \Lambda' \in P_{+}.$$
Thus, we have the following isomorphisms
 of $k$-vector spaces:
\[A_{+} \simeq \oplus_{\Lambda \in P_{+}}L(\Lambda)^{\sharp} \]
\[A_{-} \simeq \oplus_{\Lambda^{\prime} \in
P_{+}}L(\Lambda^{\prime})^{\circ}.\]
 \end{rem}

We are ready to prove the important ``triangular-type'' structure Theorem for
 $A$:

 \begin{thm}
 The multiplication map $\Delta^{*}\mid_{A_{+} \otimes A_{-}}:A_{+} \otimes A_{-} \longrightarrow
 A$
 is an isomorphism  of $U-U$-bimodules.
 \end{thm}

 \begin{pf}
Using the $U-U$ -bimodule structure on $U^{*}$ and the fact that $\Delta$ is an
algebra homomorphism on $U$, it is easily seen that $\Delta^{*}:U^{*} \otimes U^{*} \longrightarrow U^{*}$ is a $U-U$ -bimodule morphism.  Thus it suffices to prove that  $\Delta^{*}\mid_{A_{+} \otimes A_{-}}$ is bijective.
 Fix $\Lambda,\Lambda^{\prime} \in P_{+}$. Let
  $m = \Delta^{*} \mid _{L^{+}(\Lambda)^{*} \otimes  L^{-}(\Lambda^{\prime})^{*}} $.

 (1) $m$ is injective.\\
 Let $u_{\Lambda}$ and $u_{{\Lambda}^{\prime}}$ be the highest weight
 vectors of
 $L(\Lambda)$ and $L(\Lambda^{\prime})$ respectively and let
 $\{\xi_{i}\}$, $\{\xi_{j}^{\prime} \}$ be bases for $L(\Lambda)^{\sharp}$,
 $L(\Lambda^{\prime})^{\sharp}$ respectively.\\
 It suffices to prove that the
 $\{
 C_{\xi_{i},u_{\Lambda}}^{\Lambda}(C_{\xi_{j}',u_{{\Lambda}^{\prime}}}^{\Lambda^{\prime}})^{\star} \} $
 are linearly independent, which easily follows from the fact that
$$ C_{\xi_{i},u_{\Lambda}}^{\Lambda}(C_{\xi_{j}',u_{{\Lambda}^{\prime}}}^{\Lambda^{\prime}})^{\star}(U)=(C_{\xi_{i},u_{\Lambda}}^{\Lambda}{\otimes}(C_{\xi_{j}',u_{{\Lambda}^{\prime}}}^{\Lambda^{\prime}})^{\star})(\Delta (U))=
\left< \xi_{i} \otimes \xi_{j}',(1 \otimes \omega)(\Delta(U))v_{{\Lambda}}\otimes v_{\Lambda}' \right> $$
 $$= (\xi_{i} \otimes \xi_{j}')(L(\Lambda) \otimes
 L(\Lambda')),$$
since $$((1 \otimes \omega)(U))(u_{{\Lambda}}\otimes u_{\Lambda '})
=((1 \otimes \omega)(U^{\leq 0}U^{\geq 0}))(v_{\Lambda}\otimes v_{{\Lambda}'})
=((1 \otimes w)(U^{\leq 0})(1 \otimes \omega)(U^{\geq 0}))(v_{\Lambda}\otimes v_{{\Lambda}'})$$
$$=(1 \otimes w)(U^{\leq 0})((1 \otimes \omega)(U^{\geq 0})(v_{\Lambda}\otimes v_{{\Lambda}'}))
=((1 \otimes w)(U^{\leq 0}))(v_{\Lambda}\otimes L({{\Lambda}'}))
=L(\Lambda) \otimes L({\Lambda}').$$

 (2) Now we show surjectivity of the map $\Delta^{*}\mid_{A_{+} \otimes A_{-}}$.
The product of two matrix elements is given by:\[(C_{-\lambda,i;
\Lambda}^{\Lambda}(C_{-\mu,j;\Lambda^{\prime}}^{\Lambda^{\prime}})^{\star})(x)=\left<\xi_{\lambda,i}
\otimes \xi_{\mu,j}'
 ,(1 \otimes \omega)(x)(v_{\Lambda} \otimes v_{{\Lambda}'})\right>, \]
 where $\xi_{\lambda,i}$ (resp. $\xi_{\mu,j}^{\prime}$) is the
 $i$th
 (resp. $j$th) basis vector in $L(\Lambda)^{\sharp}_{\lambda}$ (resp.
 $L(\Lambda^{\prime})^{\sharp}_{\mu}$) and $x \in U$.  Recall the Remark following the definitions of $L^{+}(\Lambda)^{*}$ and  $L^{-}(\Lambda)^{*}$ in this Section.  Then, in order to prove that $\Delta^{*}\mid_{A_{+} \otimes A_{-}}:A_{+} \otimes A_{-} \longrightarrow A$ is surjective, it suffices, as in the proof of [Theorem 2.2.1] in
\cite {Yan2}, to show that the linear map:
 $$\Psi:\oplus_{\Lambda-\Lambda^{\prime}=\gamma}{\thinspace}\hom _{U}(L(\Lambda)
 \otimes
 L(\Lambda^{\prime})^{\sharp},L(\beta)) \longrightarrow
 L(\beta)_{\gamma},$$ given by
  $$\Psi(f)=f(v_{\Lambda} \otimes \xi_{-\Lambda^{\prime}})$$ is
 surjective for any $\beta \in P_{+}$ and $\gamma \in \Omega (\beta)$.
 The proof of the following result appears in \cite {Jo2}.
 \begin{lem}
 Denote by $\Psi_{\Lambda,\Lambda^{\prime}}$ the restriction of $\Psi$
 to the subspace $\hom_{U}(L(\Lambda) \otimes
 L(\Lambda^{\prime})^{\sharp},L(\beta))$. Then,
\p
  (i) $\Psi_{\Lambda,\Lambda{\prime}}$ is an isomorphism on the subspace
 $S$ of $L(\beta)_{\gamma}$, $\gamma=\Lambda-\Lambda^{\prime}$,
 consisting of the vectors $v$ such that
 $(E_i)^{n}(v)=0$ for any $n > (\Lambda,\alpha_{i}^{\vee})$ and
 $(F_i)^{n}(v)=0$ for any $n > (\Lambda^{\prime},\alpha_{i}^{\vee}) $,
 $i \in
 I$.
\p
  (ii) for $\Lambda \in P_{+}$, we have that $L(\Lambda) \simeq
 U^{\leq 0}/(U^{\leq 0})(F_i)^{(\Lambda,\alpha_{i}^{\vee})+1}$ as $U^{\leq 0}$-modules and
 $L(\Lambda^{\prime})^{\sharp}  \simeq
 U^{\geq 0}/(U^{\geq 0})(E_i)^{(\Lambda^{\prime},\alpha_{i}^{\vee})+1}$ as $U^{\geq 0}$-modules.
 \end{lem}

Since $L({\beta})$ is integrable and $ \dim L({\beta})_{\gamma} < \infty $, taking $\Lambda,\Lambda^{\prime}$ such that $\Lambda - \Lambda^{\prime}
 = \gamma$, and both $(\Lambda,\alpha_i^{\vee})$ and
 $(\Lambda^{\prime},\alpha_{i}^{\vee})$ large enough that we can get the
 whole space
 $L(\beta)_{\gamma}$ as the image of $\Psi_{\Lambda,\Lambda^{\prime}}$,
 the
 Theorem is proved.
 \end{pf}
 \begin{rem}
It easily follows that the multiplication map from $A_{-} \otimes A_{+}$ to $A$ is also bijective.
\end{rem}
 \begin{lem}
 $A$ is a domain.
 \end{lem}

 \begin{pf}
 Let $fg=0$ for some non-zero elements $f,g \in A$. Assume that $f$ and
 $g$
 are weight vectors under the $U^{0}$-$U^{0}$ action. Since the action of
$U^{>0}$-$U^{<0}$ on $U^{0}$-$U^{0}$ weight vectors in $A$ is by locally
nilpotent skew derivations, one may assume that $f$ and $g$ are $U^{\geq 0}$-$U^{\leq 0}$ invariant.
 Then there exist $\Lambda,\Lambda' \in P_{+}$ such that
 $f=C_{-\Lambda,\Lambda}^{\Lambda}$ and
 $g=C_{-\Lambda',\Lambda'}^{\Lambda'}$.  Consequently, $(fg)(1)=f(1)g(1)$, which gives the required contradiction.
 \end{pf}

 \subsection{Commutativity Relations}.
 Define $V^{>0}={\ker}(\epsilon |_{U^{>0}})$ and
 $V^{<0}={\ker}(\epsilon|_{U^{<0}})$.  For $a_{\mu} \in U_{\mu}^{>0}$,
 $b_{-\nu} \in U_{-\nu}^{<0}$,
 the following formulae are easily checked:
$$\Delta(a_\mu)=a_{\mu}{\otimes}1+K_{\mu}{\otimes}a_{\mu}{\MOD}V^{>0}\otimes
 V^{>0}$$
$$\Delta(b_{-\nu})=1{\otimes}b_{-\nu}+b_{-\nu}{\otimes}K_{-\nu}{\MOD}V^{<0}{\otimes}V^{<0}.$$

 \begin{lem} The commutation relations between elements of $A_{+}$ and $A_{-}$
are given by:
 \p
 (i) $\forall \Lambda$, $\Lambda' \in P_+$, there exist constants
 $a_{\gamma}$
 such that
 \[(C_{-\lambda,\Lambda'}^{\Lambda'})^{\star}(C_{-\mu,\nu}^{\Lambda})=q^{(\nu,\Lambda')-(\lambda,\mu)}(C_{-\mu,\nu}^{\Lambda})(C_{-\lambda,\Lambda'}^{\Lambda'})^{\star}+\sum_{\gamma}
 a_{\gamma}(C_{l_{\gamma},v_{\Lambda}}^{\Lambda})(C_{l_{\gamma'},v_{\Lambda'}}^{\Lambda'})^{\star}\]
 for $l_{\gamma} \in (l_{\mu}U^{+})_{\mu-\gamma}$ and $l_{\gamma'} \in
 (l_{\lambda}U^{+})_{\lambda-\gamma'}$.\\
 \p
 $(ii)$ Let $J_{\Lambda}(\mu,\nu)$ be the smallest 2-sided ideal of
 $A$
 containing the elements\\
 $\{ C_{-\gamma,i;\nu} ^{\Lambda} \}_{\gamma \in \Omega(\Lambda),1 \leq i
 \leq dim(L(\Lambda)_{\gamma})} $ such that $\gamma < \mu$.   Then the
 following relation holds in $A/J_{\Lambda}(\mu,\nu)$, for any $\mu,\nu \in
 \Omega(\Lambda)$ and $\lambda \in \Omega(\Lambda')$:
$$C_{-\lambda,\Lambda'}^{\Lambda'}C_{-\mu,\nu}^{\Lambda}=q^{(\nu,\Lambda')-(\lambda,\mu)}C_{-\mu,\nu}^{\Lambda}C_{-\lambda,\Lambda'}^{\Lambda'},$$
 where we use the same symbols for elements of $A$ as for their images
 in $A/J_{\Lambda}(\mu.\nu)$ under the canonical projection map.
 \end{lem}
 \begin{pf}
 $(i)$ Let $\Delta_{+}$ denote the set of positive roots of the Lie algebra
$\fg=\fg(C)$.  The following expression for the universal
 quasi-$R$-matrix for $U_{q}(\fg)$ can be found in \cite {KT}: $$R=\prod_{\alpha \in
 \Delta_{+}}{\exp}_{q_{\alpha}^{-2}}(C_{\alpha}(q)E_{\alpha}{\otimes}F_{\alpha})q^{t_{0}},$$
 where $t_{0}=\sum_{i,j}d_{ij}h_{i}{\otimes}h_{j}$ and $((d_{ij}))$ is an inverse matrix for the
symmetrical Cartan matrix $C$ if $C$ is not degenerate.  In the
case of degenerated $C$, we extend it to a non-degenerated matrix
and then take an inverse to this extended matrix. Here
$C_{\alpha}(q)$ are certain constants and  the $q$-exponential is
defined as $${\exp}_{q}(x)=\sum
_{k=0}^{\infty}q^{k(k+1)/2}x^{k}/[k]_{q}! \thinspace \thinspace
.$$ It has been established,
\cite {Gav} that in the case of irreducible integrable highest
weight $U$-modules $L$ and $L'$, the term $q^{t_{0}}$  is the
operator which acts as the scalar $q^{(\lambda,\mu)}$ on the
subspace $L_{\lambda} \otimes L'_{\mu}$.  It follows that $(1
\otimes \omega)q^{t_{0}}$ acts as the scalar $q^{(\lambda,-\mu)}$,
while
 $(1 \otimes \omega)q^{-t_{0}}$ acts as the
scalar $q^{-(\lambda,-\mu)}$, on the subspace $L_{\lambda} \otimes L'_{\mu}$.
\p
 For any $x \in U$, one has:
 $$(C_{-\lambda,\Lambda'}^{\Lambda'})^{\star}(C_{-\mu,\nu}^{\Lambda})(x)=\left<(l_{\lambda}{\otimes} l_{\mu}),(w {\otimes} 1)
 \Delta(x)v_{\Lambda'}{\otimes}v_{\nu}\right> $$
 $$=\left<(l_{\mu}{\otimes} l_{\lambda}),(1 {\otimes}
 \omega)\Delta'(x)v_{\nu}{\otimes}v_{\Lambda'}\right>
 =\left<l_{\mu}{\otimes} l_{\lambda},(1 {\otimes} \omega)R \Delta (x)
 R^{-1}v_{\nu}{\otimes}v_{\Lambda'}\right> \hspace{5.5cm}$$
$$=q^{(\Lambda',\nu)-(\lambda,\mu)}\left<(l_{\mu} {\otimes}
 l_{\lambda}),(1 {\otimes} \omega) \Delta(x) (v_{\nu} \otimes
 v_{\Lambda'})\right> + \sum_{\gamma}a_{\gamma}\left<(l_{\gamma}
 \otimes
 l_{\gamma'}),(1  \otimes \omega)\Delta (x) (v_{\nu} \otimes
 v_{\Lambda'})\right>,$$ as required.\\
 $(ii)$ is proved in a similar way.
 \end{pf}
 \begin{rem}
Note that the remaining commutativity relations in $A$ can be
derived easily, by using the involution $\star$ on the relations
(i) and (ii) above, along with Theorem 3.3.
\end{rem}

\begin{defn}
Let $P_{0}= \{ \Lambda \in P_{+}$, such that $\Lambda (h_i)=0$ for
$i \in I \}$,  $P_{\perp} = P_{+} \setminus P_{0}$, $A_{0} =
\oplus _{\Lambda \in P_{0}} C_{-\Lambda,\Lambda}^{\Lambda}$ and
$A_{\perp}= \oplus _{\Lambda \in P_{\perp}}
C_{-\Lambda,\Lambda}^{\Lambda}$.
\end{defn}

We shall henceforth assume the following: \\

{\bf Assumption A0} {\em The integrable $U$-module $L(\Lambda)$ is finite dimensional $\iff \Lambda \in P_{0}$}

\p
This is equivalent to requiring that every connected component of the Dynkin diagram is of infinite type.  Recall that the level of a highest weight module is constant.  Then the Lemma implies that the elements in $A_{\perp}$  are central in $A$.  The proof of the next proposition is the same as for Lemma 6.2.1 in \cite {Kashi1}.

\begin{prop}
The subspace $A_{\perp}$ is a 2-sided ideal of $A$.
\end{prop}

Hence we have a surjection $f:A \longrightarrow A/A_{\perp} \simeq A_{0}$.  Then, $\epsilon _{A} \circ f$ gives a one-dimensional $A$-module, denoted $N_{\infty}$, with kernel $A_{\perp}$.
 \subsection{Filtration and Gradings on $A$}
\p In the case where the Cartan matrix $C$ is of finite type, the
quantum function algebra ${\Bbb C}_q[G]$ can be filtered (see
1.4.8, 9.2.4, 10.1.4 in \cite {Jo1}).  We now construct a similar
filtration in the infinite type case. \p Recall (Section 3.3) that
$A_{-} \simeq \oplus_{\Lambda^{\prime} \in
P_{+}}L(\Lambda^{\prime})^{\circ}$. Consider $A_{-}$ as a right
$U$-module via the antipodal map $S$ on $U$ (see Section 2.3) and
write $(D_{-\Lambda,\Lambda}^{\Lambda})^{\star}$ for
$(C_{-\Lambda,\Lambda}^{\Lambda})^{\star}$.
 Let ${\cal F}^{i}(A_{-})$ be the subspaces of $A_{-}$ defined inductively
 via:
 \[{\cal F}^{0}(A_{-})=\oplus_{\Lambda \in
 P_+}{k(D_{-\Lambda,\Lambda}^{\Lambda})^{\star}}.\]
 \[{\cal F}^{i}(A_{-})={\cal F}^{i-1}(A_{-})+\sum_{\alpha \in
 \pi(C)}{\cal F}^{i-1}(A_{-})E_{\alpha}.\]Then each ${\cal F}^{i}(A_{-})$ is clearly a $U^{0}-U^{0}$-stable subspace of $A_{-}$.  Recall the fact that $\Delta (F_i)=F_i \otimes K_i^{-1} + 1 \otimes
 F_i$ and
 $(D_{-\Lambda,\Lambda}^{\Lambda})^{\star}U^{+}=L^{-}(\Lambda)^{*}$ (since $l_{-\Lambda}U^{+}=L({\Lambda})^{\circ}$), for any $\Lambda \in P_{+}$.  It follows that  ${\cal F}^{i}(A_{-})$ is a ${\Bbb Z}_{\geq 0}$-filtration of $A_{-}$.  Let us now show that $({\cal F}^{i}A_{-})A_{+}$ forms a
 filtration of $A$.
\p Let us fix $\Lambda, \Lambda' \in P_{+}$.  Then $3.4(ii)$
gives, by using induction on $i$, that for any
$C_{-\xi',\Lambda'}^{\Lambda'}$ in $A_{+}$, and
$(D_{-\xi,\Lambda}^{\Lambda})^{\star} \in {\cal F}^{i}(A_{-})$,
\[(D_{-\xi,\Lambda}^{\Lambda})^{\star}C_{-\xi',\Lambda'}^{\Lambda'}=q^{(\Lambda^{\prime},\Lambda)-(\xi',\Lambda)}C_{-\xi',\Lambda'}^{\Lambda}(D_{-\xi,\Lambda}^{\Lambda})^{\star}{\MOD}({\
 \cal F}^{i}A_{-})A_{+}.\]

Thus, in particular, $$( {\cal F}^{i}(A_{-}))A_{+}=A_{+}({\cal F}^{i}(A_{-})) \MOD ({\cal F}^{i-1}(A_{-}))A_{+}.$$

\p
We conclude,  by Theorem (3.3), that these subspaces form a filtration of $A$.

\subsection{The Prime Spectrum of $A$}
  Let $P$ be an ideal in $A$. We say that $P$ is {\em prime} in $A$ if the following condition holds: for any ideals $I_{1},I_{2}$ in $A$, $I_{1}I_{2} \subset P$ implies $I_{1} \subset P$ or $I_{2} \subset P$.
 \begin{defn}
 For $\Lambda \in P_+$, define
 \[C_{P}^{+}(\Lambda) = \{- \nu \in \Omega(\Lambda) \mid
 C_{-\nu,\Lambda}^{\Lambda} \notin P \}\]
 \[C_{P}^{-}(\Lambda) = \{- \nu \in \Omega(\Lambda) \mid
 (C_{-\nu,\Lambda}^{\Lambda})^{\star} \notin P \}.\]
If $C_{P}^{+}(\Lambda)=\emptyset$, set $D_{P}^{+}(\Lambda) = \emptyset$.
 If $C_{P}^{+}(\Lambda) \neq \emptyset$ let $D_{P}^{+}(\Lambda)$ denote
 the
 set of maximal elements of $C_{P}^{+}(\Lambda)$; define
 $D_{P}^{-}(\Lambda)$ similarly.  Fix $\vec w$=$(w_+,w_{-}) \in W \times W$.
 Define $B(w_{+},w_{-})$ as the set of all prime ideals $P$ in $A$ for
 which
 $D_{P}^{+}(\Lambda)$=${w_{+}\Lambda}$ and
 $D_{P}^{-}(\Lambda)=w_{-}\Lambda$,
 $\forall \Lambda \in P_{+}$.
 \end{defn}
 \begin{defn}
 Let $J_{\Lambda}^{+}(\lambda,\mu)$ be the 2-sided ideal in $A$
 generated by elements $C_{-\lambda',\mu'}^{\Lambda}$ with either
 $\lambda'
 \leq \lambda$ and $\mu' > \mu$ or $\lambda' < \lambda$ and $\mu' \geq
 \mu$, while $J_{\Lambda}^{-}(\lambda,\mu)=(J_{\Lambda}^{+}(\lambda,\mu))^{\star} \subset A$.
 \end{defn}
{\bf Assumption (A1)}: The collection $\Omega$ of primes $P$ in $A$ for which $D_{P}^{+}(\Lambda_{i}) \neq \emptyset$ is nonempty.
 Let Spec$_{\Omega}(A)$ denote the part of the prime spectrum of $A$ that lies in $\Omega$.
 \begin{lem}
 Assuming (A1) holds, one has
 Spec$_{\Omega}(A)=  \coprod_{(w_{+},w_{-}) \in W \times W} B(w_{+},w_{-}).$
 \end{lem}

 \begin{pf}
 Assuming that $D_{P}^{+}(\Lambda_{i}) \neq \emptyset$, pick $\eta_{i} \in
 D_{P}^{+}(\Lambda_{i})$ for $i=1,2$.
 It follows from the definition of $D_{P}^{+}(\Lambda_{i})$ that $C_{-\eta_{i},\Lambda_{i}}^{\Lambda_{i}} \notin P$ and $C_{-\nu_{i},\Lambda_{i}}^{\Lambda_{i}} \in P$, $ \forall \nu_{i} < \eta_{i}$. Thus, $J_{\Lambda_{i}}^{+}(\eta_{i},\Lambda_{i}) \subseteq P$. The commutation
 relations imply that the image $\overline{c_i}$ of
 $C_{-\eta_{i},\Lambda_{i}}^{\Lambda_{i}}$ in $A/P$ is normal, meaning
 that $$\overline{c_{i}^{\thinspace}}\thinspace\overline{c_{i}^{\star}}=\overline{c_{i}^{\star}}\thinspace\overline{c_{i}^{\thinspace}}.$$
 But
 $\overline{c_{i}} \neq 0$ by definition of $D_{P}^{+}(\Lambda_{i})$
 and,
 since $P$ is prime, $\overline{c_i}$ is a non-zero divisor.
 Now the commutation relations imply that \[ \overline{c_1}\thinspace
 \overline{c_2}=q^{2((\Lambda_1,\Lambda_2)-(\eta_1,\eta_2))}\overline{c_1}\thinspace\overline{c_2}.\]
 If we interchange the roles of $\overline{c_1}$ and $\overline{c_2}$ ,
 we get the {\em same} exponent in $q$ which must therefore be zero. In
 other words, $(\Lambda_{1},\Lambda_{2})=(\eta_{1},\eta_{2})$. The
 following
 Proposition is well-known (see Lemma A.1.17 in \cite {Jo1}):
 \begin{prop}
 Take ${\Lambda}_1,{\Lambda}_2 \in P_{+}$. Then
 $({\lambda}_1,{\lambda}_2) \leq
 ({\Lambda}_1,{\Lambda}_2)$  $\forall {\lambda}_{i} \in \Omega(\Lambda_{i})$
 and if
 ${\Lambda}_2$ is regular equality implies ${\lambda}_1=w{\Lambda}_1$
 for some $w \in W$. Then if ${\Lambda}_{1}$ is also regular,
 ${\lambda}_{2}=w{\Lambda}_{2}$.
 \end {prop}
 \noindent For some $\Lambda \in P_{++}$, it follows from the preceding
 Proposition that
 $D_{P}^{+}(\Lambda)={w\Lambda}$ for some $w \in W$. Thus $w\omega_i \in
 C_{P}^{+}(\omega_{i})$ for all $i$.
 Assume that $w\omega_{i} \notin D_{P}^{+}(\omega_{i})$ for some $i$.
 Then
 there exists $\xi \in {\Bbb N} \pi (C) \setminus \{0\}$ such that
 $w\omega_{i}-\xi
 \in D_{P}^{+}(\omega_{i})$.  Let
 $c=C_{-w\omega_{i}-\xi,\omega_{i}}^{\omega_{i}}$. Then $c \notin P$ and
 its
 image $\overline{c}$ in $A/P$ is normal and hence a non-zero divisor.
 Therefore $cc_{w\Lambda} \notin P$.  This gives that
 $w(\omega_{i}+\Lambda)-\xi \in D_{P}^{+}(\omega_{i}+\Lambda)$. Yet
 $\omega_{i}+\Lambda$ is regular so this contradicts the previous
 result.
 Suppose that there exists $\Lambda' \in P_{+}$ such that
 $D_{P}^{+}(\Lambda') \not\supseteq \{w\Lambda'\}$. Then, using the
 preceding Proposition as above, it follows that there exists
 $w' \in W$ such that $w'\Lambda' \in D_{P}^{+}(\Lambda')\setminus
 \{w\Lambda'\}$.
 Again, this forces $w\Lambda + w'\Lambda' \in C_{P}^{+}(\Lambda +
 \Lambda')$. Consequently, $w\Lambda + w'\Lambda' \geq w(\Lambda +
 \Lambda')$ and so
 $w'\Lambda' \geq w\Lambda'$.  This contradicts the fact that
 $w'\Lambda'$
 and $w\Lambda'$ are distinct minimal elements of $C_{P}^{+}(\Lambda')$
 which proves that $D_{P}^{+}(\Lambda)$=${w_{+}\Lambda}$, with
 $w=w_{+}$. The second case can be proved in a similar fashion, by the use of
the involution $\star$ on $A$.
 \end{pf}

 \par
 \vfill \pagebreak

 \section{\bf HIGHEST WEIGHT MODULES}
 \subsection{Characters}
 \begin{defn}
 Let $\chi : A_{+} \longrightarrow  k^{*}$  be an algebra
 homomorphism
 and
 $k_{\chi}$ denote the corresponding $A_{+}$-module with
 generator
 $1_{\chi}$. $\chi$ is called a {\em character} on $A_{+}$.
 \end{defn}
 Now define, for any $\Lambda \in P_{+}$ and any character $\chi$, the set
 \[ C_{\chi}^{+}(\Lambda)= \{- \nu \in \Omega(\Lambda) \slash
 \chi(C_{-\nu,\Lambda}^{\Lambda}) \neq 0 \}.\]
 If $C_{\chi}^{+}(\Lambda)=\emptyset$, or if $C_{\chi}^{+}(\Lambda)$ is
 not bounded above, set $D_{\chi}^{+}(\Lambda)=\emptyset$.
 Otherwise let $D_{\chi}^{+}(\Lambda)$ denote the set of maximal
 elements of  $C_{\chi}^{+}(\Lambda)$.  Similarly define,\[
 C_{\chi}^{-}(\Lambda)= \{- \nu \in {\Omega}(\Lambda) \slash
 \chi((C_{-\nu,\Lambda}^{\Lambda})^{\star}) \neq 0 \} \] and let
 $D_{\chi}^{-}(\Lambda)$ denote the set of maximal elements of
 $C_{\chi}^{-}(\Lambda)$.

 Consider the following assumption:

 {\bf A.2.} Let $\chi:A_{+} \longrightarrow k^{*} $ be a non-zero
 character
 such that $D_{\chi}^{+}(\Lambda) \neq \emptyset $, $\forall  \Lambda \in P_{+}$.
 \begin{lem} Assume $\chi$ satisfies A.2. If $\chi(A_{++}) \neq 0 $, then
 there exists $w \in W$ such that
  $C_{\chi}^{+}(\Lambda)=w \Lambda$.
  Further, $w$ is independent of $\Lambda$.
 \end{lem}
 \begin{pf}
 The proof of this Lemma is similar to Lemma $(3.6)$.
 Assume $D_{\chi}^{+}(\Lambda_{1}) \neq 0$ and $C_{\chi}^{+}(\Lambda_2)
 \neq 0$ and pick $\lambda_1 \in D_{\chi}^{+}(\Lambda_{1})$ and
 $\lambda_2 \in C_{\chi}^{+}(\Lambda_{2})$. Let
 $c_i=C_{-\lambda_i,\Lambda_i}^{\Lambda_{i}}$, for $i=1,2$.\\
 Then, it follows by definition that $\chi(J_{\Lambda_{1}}(\lambda_{1},\mu_{1}))=0$, and so $3.4(ii)$ gives that
\[ \chi(c_2)\chi(c_1)=q^{(\Lambda_1,\Lambda_2)-(\lambda_1,\lambda_2)}\chi(c_1)\chi(c_2).
 \]
 But $\chi(c_i)$ is a non-zero scalar, so the exponent of $q$ must be zero.
Then, Proposition $(3.6)$ implies that
 $\lambda_i=w\Lambda_i$ for a
 unique $w \in W$.  If $\Lambda$ is regular,  $D_{\chi}^{+}(\Lambda) \neq \emptyset$ implies $ C_{\chi}^{+}(\omega_{i}) \neq \emptyset $ for all $i$ by the Cartan multiplication
rule.
 \end{pf}
\begin{rem}
 Fix $w \in W$ and let $E(\Lambda,w)$ be the $w\Lambda$ weight subspace of $L^{+}(\Lambda)^{*}$ considered as a right $U^{0}$-module.  Define $I_{w}:=\oplus_{\Lambda \in P_{+}} L^{+}(\Lambda)^{*w\Lambda}$,
where $L^{+}(\Lambda)^{*w\Lambda}$ denotes the unique $U^{0}$-stable complement
in $L^{+}(\Lambda)^{*}$ of $E(\Lambda,w)$. Using the left $U^{0}$-action gives
$ L^{+}(\Lambda)^{*w\Lambda} L^{+}(\Lambda')^{*w\Lambda'} \subset L^{+}(\Lambda+\Lambda')^{*w(\Lambda+\Lambda')}$ for all $\Lambda,\Lambda' \in P_{+}$ and so the sum $I_{w}$ is a 2-sided ideal of $R^{+}$.  Thus, each $l$-tuple $\chi_{w}=\{ {\chi}_{w,i} \}_{i \in I} \in k^{l}$ can be viewed
as the unique character on $A_{+}$ satisfying ${\ker}{\chi}_{w,i}  \supset I_{w}$ and $\chi_{w}(C_{-w{\omega}_{i},{\omega}_{i}}^{\omega_{i}})=\chi_{w,i}$ for all $i$. Further, $\chi_{w}$ does not vanish on $A_{++}$ if and only if $\chi_{w,i} \neq 0$ for all $i$.  The Lemma implies all characters with this non-vanishing property are so obtained.
\end{rem}

 \subsection{The Induced Highest Weight Modules and their Irreducible Quotients $N(w)$}
Let $k(q)v_{\chi}$ denote the one-dimensional $A_{+}$-module corresponding
to each character $\chi$ on $A_{+}$, and $V(\chi):=A \otimes _{A_{+}} k(q)v_{\chi}$ the induced $A$-module.  Here $A$ acts on itself by left
 multiplication
 and $A_{+}$ acts on $A$ by right multiplication.

\begin{defn}
 We say that an $A$-module $V$ is a {\em Highest Weight Module} with {\em Highest Weight}
 $\chi$ if there
 exists a $v \in V$ such that $fv=\chi(f)v$,$ \forall f \in A_{+}$,
 and $V=Av$.

 \end{defn}
 Write $V(w)$ for $V(\chi_{w})$, where $w$ is given by Lemma $(4.1)$ and
 $\chi(C_{-w\omega_{i},\omega_i}^{\omega_{i}})=\chi_{w,i}$, $\forall i \in
 I$. The {\em highest weight vector} is denoted $v_{w}$.
 \begin{lem}
 The  ${\Bbb Z}_{\geq 0}$ graded $A_+$-module $gr_{\cal F}V(w)$
 is a
 direct sum of one-dimensional modules with characters lying in the set
 of $l$-tuples $q^{(\omega_{i},w^{-1}\eta-\Lambda)}\chi_{w,i}$, where $\Lambda \in
 P_{+}$ and $\eta \in \Omega(L(\Lambda))$.  Further, the module with character
$\chi_{w,i}$ occurs with the multiplicity of one.
 \end{lem}
 \begin{pf}
 If $(C_{-\eta,\Lambda}^{\Lambda})^{\star} \in {\cal F}^{i}A_{-}$, then
 \[(C_{-\eta,\Lambda}^{\Lambda})^{\star}C_{-\xi,\Lambda'}^{\Lambda'}=q^{(\Lambda^{\prime},\Lambda)-(\xi,\eta)}C_{-\xi,\Lambda'}^{\Lambda}(C_{-\eta,\Lambda}^{\Lambda})^{\star}{\MOD}({\
 \cal F}^{i}A_{-})A_{+}.\]
 Let $x \in gr_{\cal F}V(w)$ be such that $x={\cal
 F}^{m-1}(A_{-}){v_{w}}+y \in {\cal F}^{m}(A_{-}){v_{w}}/ {\cal F}^{m-1}(A_{-}){v_{w}}$ for
 some $m \in {\Bbb Z}$ with $y=(C_{-\eta,\Lambda}^{\Lambda})^{\star}v_{w}$.
 Then
 $$C_{-{\lambda},{\Lambda}'}^{{\Lambda}'}y=q^{(\lambda,\eta)-({\Lambda}',{\Lambda})}(C_{-\eta,\Lambda}^{\Lambda})^{\star}C_{-{\lambda},{\Lambda}'}^{{\Lambda}'}v_{w}{\MOD}{\cal
 F}^{m-1}(A_{-})A_{+}v_w.$$
 Thus, $$C_{{-\lambda},{\Lambda}'}^{{\Lambda}'}x={\cal
 F}^{m-1}(A_{-}){v_{w}}+q^{(\lambda,\eta)-({\Lambda}',{\Lambda})}(C_{-\eta,\Lambda}^{\Lambda})^{\star}C_{{-\lambda},{\Lambda}'}^{{\Lambda}'}v_{w}$$
 Take ${\Lambda}'={\omega}_{i}$ and $\lambda=w{\omega}_{i}$, together
 with the
 fact that $gr_{\cal F}V(w)=\oplus_{m \in {\Bbb Z}}gr_{m}V(w)$
  where $gr_{m}V(w)={\cal F}^{m}(A_{-}){v_{w}}/ {\cal
 F}^{m-1}(A_{-}){v_{w}}$ to
  complete the proof for the first part of the Lemma.  For the last part,
 note that the exponent in $q$ vanishes for all $i$ implies $\eta = w\Lambda$.
Yet, the corresponding element
$(C_{w\Lambda,\Lambda}^{\Lambda})^{\star}$ acts on $v_{w}$ by a
scalar.
 \end{pf}
\par
We are now ready to prove the following important theorem.
\par
\vfill\pagebreak

 \begin{thm}  The following statements are true:
\blist
\item The induced module $V(\chi)$ is a Highest Weight Module.
\item Every highest weight module $H$ of highest weight $\chi$ is an
 image of $V(\chi)$ under a surjective $A$-module homomorphism
 $\overline{\psi}:V(\chi) \longrightarrow H$.
\item $V(w)$ has a unique maximal proper submodule $V'$ and a unique irreducible quotient $N(w)$.  $N(w)$ is the unique, simple $A$-module generated by a one-dimensional
$A_{+}$-module with character $\chi_{w}$.
\elist
 \end{thm}
 \begin{pf}
\blist
 \item Define $v_{\chi}=1\otimes1_{\chi}$. Then
 $A_{+}v_{\chi}=\chi(A_{+})v_{\chi}$ which proves (i).
 \item The proof of this statement is standard.
\item  Define $V'$ as the sum of all proper submodules of $V(w)$. In
 other words, $V'$ is the sum of all submodules of $V(w)$ which do no
 contain the highest weight vector $v_{\chi}$. The preceding Lemma implies
that no proper $A$-submodule of $V(w)$ contains a copy of the one-dimensional
$A_{+}$-module with character $\chi _{w}$.  Hence the sum of any two proper $A$-submodules is again proper.  It follows that $V'$ is the maximal {\em proper} submodule of $V(w)$.  Let $N$:=$N(w)=V(w)/V'$.
 $N$ is clearly irreducible. It unique since $V'$ is unique.  As $N$ is cyclic, the
last part of the statement is easily established.
\elist
\end{pf}

\begin{rem}
If we define $J=J_{+}+J_{-}$, where $J_{+}= \sum_{\Lambda \in
P_{+}}J_{\Lambda}(\Lambda,\Lambda)$ and $J_{-}=\star(J_{+})$, then
each $l$-tuple $\chi=\{\chi_{i}\}_{i \in I}$ can be viewed as the
character on $A$ satisfying Ker$\chi \supset J$ and
$\chi(C_{-\omega_{i},\omega_{i}}^{\omega_{i}})=\chi_{i}$.  Let
$N(\chi)$ denote the corresponding one-dimensional $A$-module.
The preceding Lemma implies that $N(\chi) \simeq N(\chi_{e,i})$
given $\chi_{i}=\chi_{e,i}$ for all $i$, where $e$ is the identity
element of the Weyl group. More generally, for $w \in W$ define
$N'(w)$ to be $N(w)$ when $\chi_{w,i}=1$ for all $i$.  Then, the
preceding Lemma implies that $$N(w) \simeq N'(w) \otimes N(\chi)
\simeq N(\chi) \otimes N'(w).$$
\end{rem}
\subsection{Annihilators of $N(w)$}
 Our next goal is to describe the annihilator $J(w)$ of $N(w)$ in $A$.
 A (left) primitive ideal of $A$ is the annihilator of a
(left) simple $A$-module.
 Let Prim($A$) denote the collection of all primitive
 ideals in $A$.    It is well-known ((4.4.1) in \cite {Jo1}) that the $U$-module $L(\Lambda)_{w\Lambda}$ is 1-dimensional.
 \begin{defn} For $w\in W$, $\Lambda \in P_{+}$, let $u_{w\Lambda}$
 denote the weight vector in $L(\Lambda)$ of weight $w\Lambda$. Recall
 the quantum
 Demazure module
 $L_{w}^{+}(\Lambda)=U^{\geq 0}u_{w\Lambda}$ and let
  $L_{w}^{+}(\Lambda)^{\perp}=\{ C_{\xi,\Lambda}^{\Lambda} \mid \xi (L_{w}^{+}(\Lambda))=0$ for $\xi \in L(\Lambda)^{\sharp} \}$ be the orthogonal complement of
$L_{w}^{+}(\Lambda)$ in
  $L(\Lambda)^{\sharp}$ identified with (a subspace of)
 $L^{+}(\Lambda)^{*}$.  Similarly,
 let  $L_{w}^{-}(\Lambda)^{\perp}=\star (L_{w}^{+}(\Lambda)^{\perp})
 \subseteq L^{-}(\Lambda)^{*}.$  Now define
 $$ Q_{w}^{+}=\sum_{\Lambda \in P_{+}}L_{w}^{+}(\Lambda)^{\perp},$$
 $$ Q_{w}^{-}=\sum_{\Lambda \in P_{+}}L_{w}^{-}(\Lambda)^{\perp}.$$
 For $\vec w = (w_{+},w_{-}) \in W \times W$, let
 $Q_{(w_{+},w_{-})}$ be the minimal 2-sided ideal of $A$ containing
 $Q_{w_{+}}^{+}A_{-}$ and $A_{+}Q_{w_{-}}^{-}$. Finally, define
 $$J_{w}^{+}=\sum_{\Lambda \in P_{+}}J_{\Lambda}^{+}(w\Lambda,\Lambda)$$
 and
 $$J_{w}^{-}=\star(J_{w}^{+}).$$
 \end{defn}
\p

For any set $S \subset A$, let $\ll S \gg$ denotes the ideal in
$A$ generated by the elements in $S$.  Denote $Q_{(w,w)}$ by
$Q_{w}$ and note that
$$J_{\Lambda}^{+}(w\Lambda,\Lambda) \hspace{0.5 cm} = \hspace{0.5
cm} \ll C_{{\xi},\Lambda}^{\Lambda} \mid \xi \in
L(\Lambda)_{\lambda}^{\sharp}, \lambda \lneqq w\Lambda \gg
\hspace{0.5 cm} \subseteq \hspace{0.5 cm}\ll
L_{w}^{+}(\Lambda)^{\perp} \gg.$$ This implies that, as ideals in
$A$,
 $\ll Q_{w,w}^{+} \gg \hspace{0.2 cm} \supseteq \hspace{0.2 cm} J_{w}^{+}$.
 Let $J(w)={\ann}_{A}(N(w))$ and write the elements
$C_{-w\Lambda,\Lambda}^{\Lambda}$ as $C_{w\Lambda}$ for any $w$ and $\Lambda$.
 \begin{lem}{\bf (A)}
 $J(w) \in B(w,w)$.
 \end{lem}
 \begin{pf}
 Since $J(w)$ is primitive by definition, and primitive ideals are prime (see 3.1.5 in \cite {Dix}), it follows from Lemma $3.6$ using assumption $A1$, that we can assume $J(w) \in B(w_{+},w_{-})$.
 It suffices to prove that $w_{+}=w_{-}=w$.
 Writing $C_{-w_{+}\Lambda,\Lambda}^{\Lambda}$ as $C_{w_{+}\Lambda}$,
 relation 3.4(ii) implies that $\forall \Lambda \in P_+$, the image
 $\overline {C_{w_{+}\Lambda}} \in A \slash J(w)$ of $C_{w_{+}\Lambda}$
 under the canonical projection map is normal and nonzero.
 Thus $0 \neq \overline {C_{w_{+}\Lambda}}(A \slash J(w))v_w=(A\slash
 J(w))C_{w_{+}\Lambda}v_w$.
 Then, by Lemma $4.1$  (uniqueness of maximal element), this forces
 $w_+=w$.
 Now we prove the second case.  By definition of $w_{-}$, one has that $J(w) \supset
 \star(J_{\Lambda}^{+}(w_{-}\Lambda,\Lambda))$, for any $\Lambda \in P_+$. Let
 $d_{w\Lambda}=(C_{-w\Lambda,\Lambda}^{\Lambda})^{\star}.$  Then $3.4(ii)$
 implies
 $J_{w}^{-}(d_{w_{-}\Lambda}v_{w})=0$ and
 $C_{w\Lambda'}d_{w_{-}\Lambda}v_{w}=q^{(\Lambda',w^{-1}w_{-}\Lambda-\Lambda)}d_{w_{-}\Lambda}v_{w}$.
 If  $(\Lambda',w^{-1}w_{-}\Lambda-\Lambda)= 0$, then $w=w_{-}$, by Proposition
 3.6, and we're
 done.
 Otherwise, assume, if possible that
 $(\Lambda',w^{-1}w_{-}\Lambda-\Lambda) \neq 0$. Lemma $4.2$ implies that
 $Ad_{w_{-}\Lambda}v_{w}=A_{-}d_{w_{-}\Lambda}v_{w}$ admits no copy of
 the $A_{+}$-module $k v_{w}$ and hence is $0$ in $N(w)$. That is,
 $d_{w_{-}\Lambda}v_{w}=0$.
 Yet by $3.4$, the image of $d_{w_{-}\Lambda}$ in $A\slash J(w)$ is
 normal and non-zero. This gives the required contradiction and the
 Lemma is
 proved.
 \end{pf}
 \begin{lem}{\bf (B)}
 Every $P \in B(w_{+},w_{-})$ contains $Q_{(w_{+},w_{-})}$.
 \end{lem}
 \begin{pf}
 Given $\Lambda \in P_{+}$, $\xi \in L(\Lambda)^{\sharp}_{\eta}$ such
 that $C_{\xi,\Lambda}^{\Lambda} \notin P$ and $C_{\xi \cdot
 a,\Lambda}^{\Lambda} \in P$, $\forall a \in U_{+}^{+}$, then claim
 $\eta =
 w_{+}\Lambda$.
 Indeed, if $\nu \neq w_{+}\nu$, then there exists $a_{\zeta} \in
 {\ker} ({\epsilon}{\mid}_{U^{\geq 0}})$ of weight $\zeta=\nu-w_{+}\nu$ such
 that $\xi_{\nu} \cdot
 a_{\zeta}=\xi_{w_{+}\nu}$.
 But
 \[C_{-\nu,\nu}^{\nu}C_{-\eta,\Lambda}^{\Lambda}=q^{-(\nu,\eta)+(\nu,\Lambda)}C_{-\eta,\Lambda}^{\Lambda}
 C_{-\nu,\nu}^{\nu}.\]
 Acting on the right by $a_{\zeta}$ on both sides of this equation gives
 that
 \[C_{-w_{+}\nu,\nu}^{\nu}C_{-\eta,\Lambda}^{\Lambda}=q^{-(w_{+}\nu,\eta)+(\nu,\Lambda)}C_{-\eta,\Lambda}^{\Lambda}C_{-w_{+}\nu,\nu}^{\nu}
 {\MOD} P.\]
 Yet $P \supset J_{w_{+}}^{+}$, so the commutation relations (3.4) give
 that
 \[C_{-\eta,\Lambda}^{\Lambda}C_{-w_{+}\nu,\nu}^{\nu}=q^{(\Lambda,\nu)-(\eta,w_{+}\nu)}C_{-w_{+}\nu,\nu}^{\nu}\C_{-\eta,\Lambda}^{\Lambda}{\MOD}P.\]
 Now $C_{-w_{+}\nu,\nu}^{\nu} \notin P$ implies that it has a normal
 image which is a non-zero divisor. Thus $(w_{+}\nu,\eta)=(\Lambda,\nu)$
 and
 therefore $w_{+}\eta=\Lambda$ and the claim is proved.
 Similarly, it can be shown that given $\Lambda \in P_{+}$, $\xi \in
 L(\Lambda)^{\sharp}_{\eta}$ such that $(C_{\xi,\Lambda}^{\Lambda})^{\star}
 \notin P$ and $(C_{\xi \cdot a,\Lambda}^{\Lambda})^{\star} \in P$, $\
 \forall a \in U_{+}^{+}$, then  $\eta = w_{-}\Lambda$.
 Now suppose, on the contrary, that $P \not\supset Q_{(w_{+},w_{-})}$.
 Assume that $C_{\xi,\Lambda}^{\Lambda} \notin P$ and yet
 $C_{\xi,\Lambda}^{\Lambda} \in P$, with
 $\xi \in L(\Lambda)_{\eta}^{\sharp}$ and $\eta$ maximal with respect to this property. Since
 $Q_{w_{+}}^{+}$ is stable for the right action of $U^{+}$, it follows
 that
 $C_{\xi \cdot a,\Lambda}^{\Lambda} \in P$ for all $a \in U_{+}^{+}$.
 This
 forces $\eta = w_{+}\Lambda$, contradicting that
 $C_{-w_{+}\Lambda,\Lambda}^{\Lambda} \notin Q_{w_{+}}^{+}$.  The
 second case can be proved in a similar fashion.
 \end{pf}

 An immediate consequence of these Lemmas is the following important
 Theorem:
 \begin{thm}
 $J(w) \supset Q_{w}$
 \end{thm}

 \subsection{Correspondence with Schubert cells}
 Fix $\Lambda \in P_{+}, w \in W$ of length $l(w)$ and recall (Section 4.3) that the ideal generated by $Q_{w}^{+}$ in $A$ contains $J_{w}^{+} \supset  J(w,\Lambda)$
where $J(w,\Lambda)=J_{\Lambda}^{+}(w\Lambda,\Lambda)$ is
 the 2-sided
 ideal in $A$ generated by $\{ C_{-\mu,i;\Lambda}^{\Lambda}$ with  $\mu \lneq w\Lambda \}$ and $J_{w}^{+}$=$\sum_{\Lambda \in P_{+}}J(w,\Lambda)$.
 As in \cite {KP}, a {\em Schubert Cell} is defined as $S_{w}=B_{+}wB_{+}/B_{+}$ in $G/B_{+}$, with dim$S_{w}=2l(w)$.
 \begin{defn}
 Let $\phi:A\longrightarrow {\End}(V)$ be an irreducible representation of
 $A$.
 We say that $(V,\phi)$ {\em corresponds to the Schubert Cell} {\bf
 $S_{w}$} if the following two conditions are satisfied:\\
 1. $\phi(Q_{w}^{+})=0$, $\forall \Lambda \in P_{+}$.\\
 2. $\phi(C_{-w\Lambda,\Lambda}^{\Lambda}) \neq 0$,  $ \forall \Lambda
 \in P_{+}$.
 \end{defn}

 \begin{rems}
\blist
\item Some reasons for this terminology are discussed now.  Classically, a {\em Flag Manifold} $F_{\Lambda}$ is the orbit in the
 projective space ${\Bbb P}(L(\Lambda))$ of a highest weight vector.
 Given
 $\Lambda \in P_{+}$, the {\em quantized algebra of functions on the
 Flag
 Manifold $F_{\Lambda}$}, also called the Representation ring,  is the
 subalgebra of $A$ generated by the matrix coefficients of the form $ \{
 C_{-\mu,i;\Lambda}^{\Lambda} , \forall \mu,i \}$. Denote it as
 $A[F_{\Lambda}]$. While observing that
 it is not a $\star$ -subalgebra of $A$, its algebra structure can be described
 explicitly via the Cartan multiplication rule (projection on highest
 weight component of the tensor product $L(m\Lambda)^{\sharp} \otimes
 L(n\Lambda)^{\sharp} \longrightarrow L((m+n)\Lambda)^{\sharp}$, since
 $A[F_{\Lambda}] \simeq \bigoplus _{k=0}^{\infty}L(k\Lambda)^{\sharp}$
 as $U$-modules). Classically, (see Section 2.6 in \cite {KP}) the {\em Schubert varieties} $S_{w}^{\Lambda}$ are
 the closures in the projective space ${\Bbb P}(L(\Lambda))$ of the
 $B_{+}$-
 orbits of the extremal vector of weight $w\Lambda$, which is unique
 up to a scalar multiple. If $\Lambda \in P_{++}$ is strictly dominant,
 then $S_{w}^{\Lambda} \simeq S_{w}$ is the closure of a Schubert cell.
 Now, for $\Lambda \in P_{+}$,
 $w \in W$, let $\dot Q(w,\Lambda)$ be the 2-sided ideal in
 $A[F_{\Lambda}]$ generated by the matrix coefficients
 $C_{-\mu,i;\Lambda}^{\Lambda}$
 such that $\mu \ngeq w\Lambda$. Then the
 quotient algebra {\bf $A[F_{\Lambda}]\slash\dot Q(w,\Lambda)$}
 is called the {\em quantized algebra of functions on the Schubert
 variety} $S_{w}^{\Lambda}$ and is denoted by  $A[S_{w}^{\Lambda}]$.
 Note
 that $\dot Q(w,\Lambda) \subseteq  L_{w}^{+}(\Lambda)^{\perp}$, with equality in the $sl_2$ case.
 The motivation for the definition of an irreducible
 representation corresponding to a Schubert cell is then apparent from Theorem 4.3, since such a representation induces a representation of
 $A[S_{w}^{\Lambda}]$ for all $\Lambda \in P_{+}$ and may thus be viewed as
being  ``supported'' on these quantum Schubert varieties.
\item In the affine case, the module $N_{\infty}$ (Proposition 3.4) does not correspond to any Schubert cell, since its annihilator contains all of $A_{\perp}$.  At present, we know of no other irreducible representations other than
these. It would be nice if one could establish that these modules
do indeed exhaust all the irreducible representations of
$A$.\elist
\end{rems}
\par
\vfill \pagebreak

\section{\bf WEIGHT SPACE DECOMPOSITIONS}
We will assume that $0<q<1$ for this Section (this condition is
necessary for unitarizability).
\begin{defn}
Recall the algebra $A=C[G]_q$ carries an involution $\star$.
 An $A$-module $V$ is said to be {\em unitarizable} if it admits a
positive-definite Hermitian form $\left<,\right>:V {\times} V
 \longrightarrow k$,
 such that $$\left<av_{1},v_{2}\right>=\left<v_{1},a^{\star}v_{2}\right>,$$ for any $a \in A,
   v_{1},v_{2} \in V$.
 \end{defn}

 \subsection{Weight Space Decompositions of N(w)}
 The goal now is to prove the Tensor Product
 Theorem for the structure of the irreducible $A$-modules corresponding to the
 Schubert cells.
 \begin{prop}{\bf (A)} \thinspace
  For any $\Lambda \in P_{+}$ the following identity holds:, \[
 \sum_{\mu,i}(C_{-\mu,i;\Lambda}^{\Lambda})^{\star}C_{-\mu,i;\Lambda}^{\Lambda}
 =
 1,\]
where the sum is taken over all weights $\mu$ of $L(\Lambda)$ and $1\leq i \leq \Dim (L(\Lambda))$.
 \end{prop}
 \begin{pf}
 Define $x_{\Lambda} \in A_{-} \otimes A_{+}$  as the element corresponding
 to the left-hand-side of the above sum. Using the formulae (Section 3.4):
$$\Delta(a_\mu)=a_{\mu}{\otimes}1+K_{\mu}{\otimes}a_{\mu}{\MOD}V^{>0}\otimes
 V^{>0}$$
$$\Delta(b_{-\nu})=1{\otimes}b_{-\nu}+b_{-\nu}{\otimes}K_{-\nu}{\MOD}V^{<0}{\otimes}V^{<0},$$
along with the triangular structure of $U$, one has that
 $\Delta^{*}(x_{\Lambda}(u))=\epsilon(u), \forall u \in U$. $x_{\Lambda}$ is therefore a bi-invariant element of $A$ and so it equals $1$.
 \end{pf}
 \p
Returning now to the irreducible $A$-modules $N(w)$ constructed in
Section 4.2, we let $B_{w}$ be the $\star$-subalgebra in
${\End}_{k}(N(w))={\hom}_{k}(N(w),N(w))$ generated by
 the operators corresponding to $\{C_{-w\Lambda,\Lambda}^{\Lambda}, \Lambda \in P_{+}\}.$  One easily
 checks that it is commutative,
due to the fact that  Ann $N(w) \supseteq Q_{w} \supseteq J_{w}^{+} + J_{w}^{-}$. It
 follows that the elements $C_{w\Lambda}$ and $(C_{w\Lambda})^{\star}$ have
 normal images $\overline{C_{w\Lambda}}$ and
 $\overline{d_{w\Lambda}}=\overline{(C_{w\Lambda})^{\star}}$ in
 $A$/Ann$N(w)$, which means that
$$\overline{C_{w\Lambda}^{\thinspace}}\thinspace\overline{d_{w\Lambda}}=\overline{d_{w\Lambda}}\thinspace\overline{C_{w\Lambda}^{\thinspace}}.$$

\begin{defn}
The {\em weight space} $N(w)_{\gamma}$, of weight $\gamma \in {\Bbb N}\pi(C)$, is defined as the subspace
of $N(w)$ in which the commuting elements $\overline{C_{w\Lambda}}$ and
 $\overline{d_{w\Lambda}}=\overline{(C_{w\Lambda})^{\star}}$ act by the scalar $q^{(\gamma,\Lambda)}$, for any $\Lambda \in P_{+}$.
\end{defn}

\begin{rem}
Make the important observation that a highest weight $A$-module
has highest weight $\gamma=0$, hence $q^{(\gamma,\Lambda)}=1$.
That is, $\chi=1$, where $\chi$ is the highest weight of the
induced module $V(\chi)$.  More generally (see Proposition(B)) the
weight $\gamma'$ of $N(w)$ corresponds to a weight $\chi'$ of
$V(\chi)$ via
$\chi'(C_{-w\Lambda,\Lambda}^{\Lambda})=q^{(\gamma',\Lambda)}$,
since
$(C_{-w\Lambda,\Lambda}^{\Lambda})^{\star}v_{\gamma'}=C_{-w\Lambda,\Lambda}^{\Lambda}v_{\gamma'}
=q^{(\gamma',\Lambda)}v_{\gamma'}$, for any $v_{\gamma'} \in
N(w)_{\gamma'}$, $\Lambda \in P_{+}$.
\end{rem}

\p
Let $\Omega(w)$ denote the collection of all $\gamma$ for which $N(w)_{\gamma} \neq 0$.
Then, for any $\Lambda,\Lambda' \in P_{+},$ and any $\eta \in \Omega(\Lambda)$, the commutation
 relations (3.4) imply that
 \[ \overline {(C_{-\eta,\Lambda}^{\Lambda})^{\star}}\thinspace
 \overline{C_{w\Lambda'}^{
 }}=q^{(\Lambda',\Lambda)-(\eta,w\Lambda')}\overline{C_{w\Lambda'}^{
 }}\thinspace\overline{(C_{-\eta,\Lambda}^{\Lambda})^{\star}}\]
 \[ \overline {(C_{-\eta,\Lambda}^{\Lambda})^{\star}}\thinspace
 \overline{d_{w\Lambda'}^{
 }}=q^{(\Lambda',\Lambda)-(\eta,w\Lambda')}\overline{d_{w\Lambda'}^{
 }}\thinspace\overline
 {(C_{-\eta,\Lambda}^{\Lambda})^{\star}}.\]

 \begin{prop}{\bf (B)}
 \noindent  Let $N=N(w)$ be a simple module corresponding to a
 Schubert cell $S_{w}$. Then $N=\bigoplus _{\gamma \in {\Omega}(w)}N_{\gamma}$.
 \end{prop}
 \begin{pf}
Recall that $N(w)=A_{-}v_{w}$. Then it holds that $\gamma \in \Omega(w)$ satisfies $\gamma \leq 0$. Indeed, since
$$C_{w\Lambda'}(C_{-\eta,\Lambda}^{\Lambda})^{\star}v_{w}= q^{(\Lambda',w^{-1}\eta-\Lambda)} (C_{-\eta,\Lambda}^{\Lambda})^{\star}v_{w},$$ and also, $${d_{w\Lambda'}}{(C_{-\eta,\Lambda}^{\Lambda})^{\star}}\thinspace v_{w}=q^{(\eta,w\Lambda')-(\Lambda',\Lambda)} {(C_{-\eta,\Lambda}^{\Lambda})^{\star}}\thinspace {d_{w\Lambda'}}v_{w}=q^{(w^{-1}\eta-\Lambda,\Lambda')} (C_{-\eta,\Lambda}^{\Lambda})^{\star} v_{w},$$ we have that $(C_{-\eta,\Lambda}^{\Lambda})^{\star}v_{w}$ belongs to the subspace $N(w)_\gamma$ with $\gamma=w^{-1}\eta-\Lambda$.  Thus,   $$\Omega(w)= \bigcup _{\Lambda \in P_{+}}(w^{-1}\Omega(L_{w}^{-}(\Lambda))-\Lambda).$$ Now $\gamma=0$ only when $\eta = w\Lambda$, in which case $(C_{-w\Lambda,\Lambda}^{\Lambda})^{\star}v_{w}=v_{w}$ and so
$N_{0}=kv_{w}$.  For any other $d \in A_{-}$, it follows as above that $dN_{\gamma} \subseteq N_{\gamma'}$, with $\gamma' < \gamma$.
 Theorem (4.3) implies that Ann$(N)
 \supseteq
 Q_{w}^{-}$, so $A_{-}$ can be replaced by $A_{-}/Q_{w}^{-}$. But,
 since $L(\Lambda)^{\sharp}/L_{w}^{-}(\Lambda)^{\perp} \simeq
 L_{w}^{-}(\Lambda)^{*}$,
 $A_{-}/Q_{w}^{-}$ can be identified with $\oplus_{\Lambda \in P_{+}}
 L_{w}^{-}(\Lambda)^{*}$
 which, in turn, are spanned by the $(C_{\xi,\Lambda}^{\Lambda})^{\star}$, $\xi
 \in L_{w}^{-}(\Lambda)^{*}$.

\p Since, $\Omega( L_{w}^{-}(\Lambda)^{*})=\Omega( L_{w}^{-}(\Lambda))$, the
 assertion follows from the formulae immediately preceding this
 Proposition.
 \end{pf}
 \begin{prop}{\bf (C)}
  $\forall \Lambda \in P_{+}$, and for any $\chi$ as in (4.2) it holds that
 $|\chi(C_{-w\Lambda,\Lambda}^{\Lambda})|=1$.
 \end{prop}
 \begin{pf}
 This is an immediate consequence of Proposition(A) and the fact that if
 $v \in N_{\chi}$, then $C_{-\mu,\Lambda}^{\Lambda}v=0, \forall \mu \in
 \Omega(\Lambda) \setminus \{ w \Lambda \}$.
 \end{pf}

 \par
 \vfill\pagebreak

 \section {\bf TENSOR PRODUCT THEOREM}
 Continue to assume that $0<q<1$ for this Section.
 \subsection{Elementary Modules}
 Let $\alpha_i$ be a simple root of the Kac-Moody algebra $\fg=\fg (C)$.
 Define $\psi_i:U_{q_i}(sl(2,{\Bbb C})) \rightarrow U_{q}({\fg}) $ as
 the canonical embedding of Hopf algebras given by
 $E \mapsto E_i , F \mapsto F_i, K \mapsto K_i $.
 Consider the restriction of the dual morphism $\psi_i^*:A \longrightarrow
 {k}_{q_i}[SL_2]$, which is surjective, since every finite dimensional $U_i=U_{q_i}(sl(2,{\Bbb C})$ module occurs in some integrable $U$-module (see 4.3.6 \cite {Jo1}).
 We see that, given an irreducible representation $\phi$ of  $
 {k}_{q_i}[SL_2]$, we get an irreducible representation $\phi_i = \phi
\circ \psi_i^* $ of $A$.  The corresponding $A$-module is called {\em
 elementary}, and is isomorphic to $N(s_{i})$, by $4.2(iii)$ (see 10.1.4 in \cite {Jo1}). It will be denoted simply as $N(i)$.

 \subsection {Proof of the Tensor Product Theorem}
 Consider the following Lemma:

\begin{lem}
 Let $i \in I, w \in W$ such that $l(s_{i}w)>l(w)$. Then \\
 $(i)$ ${\ann}_{A}(N(i) \otimes N(w)) \supseteq Q_{s_{i}w}^{+}$, for all $a \in N(i), b \in N(w)$ .\\
 $(ii)$ $C_{-s_{i}w\Lambda,\Lambda}^{\Lambda} (a \otimes b) =
 C_{-s_{i}w\Lambda,w{\Lambda}}^{\Lambda}(a) \otimes
 C_{-w\Lambda,\Lambda}^{\Lambda}(b)$.
 \end{lem}

 \begin{pf}

 $(i)$ The hypothesis  $l(s_{i}w)>l(w)$ implies that
 $L_{s_{i}w}^{+}(\Lambda)$ is
 $U_{i}$-stable, since $L_{s_{i}w}^{+}(\Lambda)=U_{i}L_{w}^{+}(\Lambda) \implies U_{i}L_{s_{i}w}^{+}(\Lambda)=L_{s_{i}w}{+}(\Lambda)$. Using the non-degenerate Shapovalev form on $L(\Lambda)$, identify
 $L_{s_{i}w}^{+}(\Lambda)^{\perp\
 }$ in $L(\Lambda)^{\sharp}$ with a $U_{i}$-stable subspace of
 $L(\Lambda)$ complementing $L_{s_{i}w}^{+}(\Lambda)$.
 Fix a basis $\{u_{j}\}_{j \in {\Bbb Z}_{\geq 0}}$ for $L(\Lambda)$
 which is compatible with the inclusion $L_{w}^{+}(\Lambda) \subseteq
 L_{s_{i}w}^{+}(\Lambda)$,
 i.e. with the choice of bases $\{u_{j}\}_{j=1}^{s}$ for
 $L_{w}^{+}(\Lambda)$,
 $\{u_{j}\}_{j=1}^{r}$ for $L_{s_{i}w}^{+}(\Lambda)$ (see \cite {Lak}) and
 $\{u_{j}\}_{j>r}$
 for $L_{s_{i}w}^{+}(\Lambda)^{\perp}$ .  Let
 $\{ l_{j} \} $ be the dual basis for $L(\Lambda)^{\sharp}$.
 Let $\xi \in L(\Lambda)^{\sharp}$, $a \in N({i})$, $b \in N(w)$. Then
 $$C_{\xi,\Lambda}^{\Lambda}(a \otimes b) = \sum_{j \in {\Bbb Z}_{\geq
 0}}
 C_{\xi,u_{j}}^{\Lambda}a \otimes C_{l_{j},\Lambda}^{\Lambda}b,$$
 where $$\Delta(C_{\xi,\Lambda}^{\Lambda})=\sum_{j \in {\Bbb Z}_{\geq
 0}}C_{\xi,u_{j}}^{\Lambda} \otimes C_{l_{j},\Lambda}^{\Lambda}.$$
\p
 Now, by Theorem 4.3, we have that  ${\ann}_{A} N(w) \supseteq Q_{w}^{+}$. It follows that
 $C_{l_{j},\Lambda}^{\Lambda}b=0$, for any $j>s$. Clearly
 $C_{\xi,u_{j}}^{\Lambda}(a)=\psi_{i}^{*}(C_{\xi,u_{j}}^{\Lambda})(a)$.
 For $\xi \in
 L_{s_{i}w}^{+}(\Lambda)^{\perp}$ and $v \in L_{s_{i}w}^{+}(\Lambda)$, we have
 that
 $C_{\xi,v}^{\Lambda}(x)=\xi(xv)=0$ if $x\in U_{i}$. So, in particular,
 $C_{\xi,u_{j}}^{\Lambda}a=0$   $\forall j \leq r$ and $\xi \in
 L_{s_{i}w}^{+}(\Lambda)^{\perp}$, which proves $(i)$.

 $(ii)$ Let $\xi \in L(\Lambda)^{\sharp}$ of weight $-s_{i}w{\Lambda}$.
 It is clear that the only non-zero contributions to the right hand side
 of the expression
 $$C_{\xi,\Lambda}^{\Lambda}(a \otimes b) = \sum_{j \in {\Bbb Z}_{\geq
 0}}
 C_{\xi,u_{j}}^{\Lambda}a \otimes C_{l_{j},\Lambda}^{\Lambda}b$$
 occur when $j \leq s$ and $C_{\xi,u_{j}}^{\Lambda}(x) \neq 0$
 for some $x \in U_{i}$.
\p
 Further, one can assume $u_{j}$ to be weight vectors with
 $u_{s}=u_{w{\Lambda}}$.
 But the $U_{i}$-module generated by $\xi$ has weights
 $\{-s_{i}w\Lambda,-s_{i}w\Lambda-\alpha_{i},...,-w\Lambda \}$. The
 definition of
 $L_{w}^{+}(\Lambda)$ then implies that the only non-zero term corresponds
 to $j=s$.
 So $u_{j}=u_{w\Lambda}$ whereas $l_{j}=-w\Lambda$,  which proves $(ii)$.\\
 \end{pf}
\par We are now ready for our main result, the Tensor Product Theorem:
 \begin{thm}
 $N:=N({i}) \otimes N(w)$ is unitarizable, irreducible, and isomorphic to
 $N(s_{i}w)$.
 \end{thm}
 \begin{pf}
 Observe that $N({i})$ is unitarizable, by the results of 5.1 and 6.1.
 The tensor product of any 2 unitarizable $A$-modules $M_1$ and $M_2$
 is also a unitarizable $A$-module under $$\left<v_{1} \otimes v_{2}, v_{3} \otimes v_{4}\right>=\left<v_{1}, v_{3}\right>_{M_{1}}\left<v_{2}, v_{4}\right>_{M_{2}}.$$ Thus, it follows using
 induction that $N$ is unitarizable. We shall prove the irreducibility using induction on the length of $w$.
\p
 For $l(w)=1$, $N(i)$ is irreducible by definition (6.1). Let
 $w=s_{i_{1}}s_{i_{2}}...s_{i_{k}}$ be a reduced expression for $w$.
 Suppose  $N(w) \simeq N({i_{1}}) \otimes ....
 N({i_{k}})$, and
 is irreducible.
\p
 Let $\{ e_k \}_{k=0}^{\infty}$ (resp. $\{f_{M}\}_{M \in {\Bbb
 Z}_{+}^{l(w)}} $) be an orthonormal basis of $N({i})$ (resp. $N(w)$).
 Thus
 $\{e_{k} \otimes f_{M}\} $ is an orthonormal basis for $N$.
 Then the arguments in the proof of Lemma 6.2(ii) imply that
 $$C_{-s_{i}w\Lambda,\Lambda}^{\Lambda} (e_{k} \otimes f_{M}) =
 \psi _{i}^{*}C_{-s_{i}w\Lambda,w{\Lambda}}^{\Lambda}(e_{k}) \otimes
 C_{-w\Lambda,\Lambda}^{\Lambda}(f_{M}).$$
\p
 It is easy to see (via the commutation relations) that Ker$(C_{-w\Lambda,\Lambda}^{\Lambda}) $ is an $A$-invariant subspace in
 $N(w)$. But $N(w)$
 is irreducible, which forces this kernel to be trivial.
 Hence, using $C_{-s_{i}w\Lambda,\Lambda}^{\Lambda} (e_{k} \otimes f_{M}) =
 \psi _{i}^{*}C_{-s_{i}w\Lambda,w{\Lambda}}^{\Lambda}(e_{k}) \otimes
 C_{-w\Lambda,\Lambda}^{\Lambda}(f_{M})$, Lemma 6.2 implies the following:
\p (1) The $\star$-algebra $B_{s_{i}w}$ in $\End(N)$ generated by
 $C_{-s_iw\Lambda,\Lambda}^{\Lambda} (\Lambda \in P_{+})$ is commutative
 and diagonalizable under the basis $\{ e_k  \otimes f_M \}$.
\p
(2) The eigenvalue $\chi _{\Lambda}^{(k,M)}$ of  $C_{-s_iw\Lambda,\Lambda}^{\Lambda}$ corresponding to $e_k \otimes f_M$
 satisfies $|\chi_{\Lambda}^{(k,M)}| < 1$ unless $k=0,M=M_{0}$, where
 $M_{0}=(0,0,...,0)$.
\p
(3) $e_0 \otimes f_{M_0}$ is the only eigenvector for any
 $C_{-s_iw\Lambda,\Lambda}^{\Lambda}$ with the  property that
 $|\chi_{\Lambda}^{(0,M_{0}))}|=1$
\p
 Also, all eigenvalues $\chi_{\Lambda}^{(k,M)}$ are nonzero.
 Now, let $W' \subsetneqq N$ be a non trivial $A$-invariant subspace, such
that $e_{0} \otimes f_{M_{0}} \notin W'$.  Thus, by (1), the action of $B_{s_iw}|_{W'}$ is diagonalizable.  For any $v \in N$, Proposition 5.2(A) implies that
 $$\left<C_{-\mu,\Lambda}^{\Lambda} v,C_{-\mu,\Lambda}^{\Lambda}v\right> \leq  \left<v,v\right>.$$
 Then, by Lemma (3.4), there exists a vector $v$ in $W'$ such that $fv=\chi(f)v$ for any $f
 \in A_+$, which can be chosen among the vectors of the basis $\{ e_{k}
 \otimes f_{M}\}$ with
 $(k,M) \neq (0,M_{0})$. If this is not the case, take the orthogonal
 complement of $W'$ which is also an $A$-module. Therefore
 $|\chi(C_{-s_iw\Lambda,\Lambda}^{\Lambda})| < 1$ for some $\Lambda
 \in P_{+}$.
\p
 Calling such a $v$ in $W'$ as
 {\em primitive} and the corresponding homomorphism $\chi:A_{+}
 \longrightarrow \Bbb C$ its {\em weight}, we introduce a natural partial
 order
 $\chi _{1} \prec \chi _{2} \hspace{0.2 cm} \Leftrightarrow |\chi _{1} (f)| \leq |\chi
 _{2} (f)|$, for any $f \in \{ C_{-s_{i}w\Lambda ,\Lambda}^{\Lambda}
 \}
 _{\Lambda \in P_{+}}$ (we can use the same notations, since the homomorphism $\chi$ has a canonical extension to $B_{s_{i}w}$ such that $\chi(( C_{-s_{i}w\Lambda ,\Lambda}^{\Lambda})^{\star})=\overline{\chi( C_{-s_{i}w\Lambda ,\Lambda}^{\Lambda})}$).  Let $\chi '$ be a maximal element in the set
 of
 primitive weights, with respect to
 the above order, with eigenvector $v'$. Then the commutation relations
 imply that any vector of the form $C_{-\lambda,j;\Lambda}^{\Lambda}v'$
 is an eigenvector for $B_{s_{i}w}$ with weight larger than
 $\chi '$, which is
 absurd.
 This forces $C_{-\lambda,i;\Lambda}^{\Lambda} v'=0$,  $ \forall \lambda
 \neq s_iw\Lambda$. Now Proposition $5.2(A)$ implies that
 $\sum_{\mu,j}(C_{-\mu,j;\Lambda}^{\Lambda})^{\star}C_{-\mu,j;\Lambda}^{\Lambda}v'=
 v'.$ But all except one of the summands in this equation disappear, and so
 $$(C_{-s_{i}w\Lambda,j;\Lambda}^{\Lambda})^{\star}C_{-s_{i}w\Lambda,j;\Lambda}^{\Lambda}v'=v'.$$
 One may assume that $\left<v',v'\right>_{N(i) \otimes N(w)}=1$. Thus it follows
 that $|\chi(C_{-s_{i}w\Lambda,\Lambda}^{\Lambda})|=1$,  $\forall
 \Lambda \in P_{+}$. This is
 possible only if $v'=e_{0} \otimes f_{M_{0}}$. But $e_{0} \otimes
 f_{M_{0}} \notin W'$ by assumption.
 Thus we have a contradiction. Hence $N$ is irreducible.
 The fact that $N \simeq N(s_{i}w)$ follows at once from the preceding
 Lemma and Theorem $4.2(iv)$.
 \end{pf}

 \begin{cor}
 The characterization $N(w) \simeq N({i_{1}}) \otimes ....
 N({i_{k}})$ is
 independent of the choice of reduced decomposition
 $s_{i_{1}}...s_{i_{k}}$ for $w$.
\p
 \end{cor}
 \par
 \vfill\pagebreak

 \section {\bf OPEN QUESTIONS}
\blist
\item Compute the prime and primitive spectra of ${\Bbb
 C}_{q}[G]$ as well as the
 spectrum of the Abelian category ${\Bbb
 C}_{q}[G]-mod$ in the sense of A.Rosenberg \cite {Ros}.
  \item Are there any  solutions to the Zamolodchikov
 tetrahedron
 equations, using the intertwiners of isomorphic, irreducible $A$-modules
 ?
\item Define the affine quantum Weyl group and relate it
 with
 the Dynamical Quantum Weyl groups of Varchenko and Etingof \cite
 {EV}.
 \item Compute the spectra of the quantized function
 algebra, $R={\Bbb C}_{q}[G/B_{+}]$ (resp. $R_{+}={\Bbb
 C}_{q}[G/N_{+}]$)  of the flag
 manifold (resp. base affine space).
 \item What are the answers to the above questions
 when you specialize to the roots of unity case?
 \elist

 \par
 \vfill\pagebreak

 \end{document}